\documentclass{amsproc} 
\usepackage{amssymb,amsmath,latexsym,epsfig,euscript,eepic,rotating,color,wrapfig} 

\definecolor{yell}{rgb}{1,1,.2} 
\newlength{\sh}
\newlength{\hs}
\newlength{\agk}
\newlength{\agku}
\newlength{\jmr}
\newlength{\jfc}
\newlength{\mum}
\newlength{\koi}
\newlength{\hwl}
\newlength{\khov}
\newlength{\bernd}
\newlength{\clo}
\newlength{\groth}
\newlength{\gd}
\newtheorem{lemma}{Lemma}[subsection] 
\newtheorem{nota}{Notation}[subsection]

\newtheorem{thm}{Theorem}[subsection]
\newtheorem{prop}{Proposition}[subsection]
\newtheorem{dfn}{Definition}[subsection]
\newtheorem{cor}{Corollary}[subsection]
\newtheorem{rem}{Remark}[subsection]	
 
\newtheorem{ex}{Example}[subsection]

\newcommand{\cN}{\mathcal{N}}
\newcommand{\cO}{\mathcal{O}}
\newcommand{\supp}{\mathrm{Supp}}
\newcommand{\init}{\mathrm{Init}}
\newcommand{\conv}{\mathrm{Conv}}

\newcommand{\thth}{^{\underline{\mathrm{th}}}}

\newcommand{\nd}{{\underline{\mathrm{nd}}}}
\newcommand{\Pro}{{\mathbb{P}}}

\newcommand{\R}{\mathbb{R}}
\newcommand{\C}{\mathbb{C}}
\newcommand{\N}{\mathbb{N}}
\newcommand{\Z}{\mathbb{Z}}
 
\newcommand{\cp}{\mathfrak{p}}

\newcommand{\area}{\mathrm{Area}}

\newcommand{\newt}{\mathrm{Newt}} 
\newcommand{\jac}{\mathrm{Jac}} 
\newcommand{\crit}{\mathrm{Crit}} 
\newcommand{\toricjac}{\mathrm{ToricJac}} 
\newcommand{\cay}{\mathrm{Cay}} 
\newcommand{\rank}{\mathrm{rank}} 
 
\newcommand{\perm}{\mathrm{Perm}}

\newcommand{\Zn}{\Z^n}

\newcommand{\Rn}{\R^n}

\newcommand{\Cn}{\C^n} 
 
\newcommand{\Cs}{\C^*}

\newcommand{\Csn}{{(\C^*)}^n}
\newcommand{\gln}{\mathbb{G}\mathbb{L}_n}
\renewcommand{\qed}{$\blacksquare$}
\newcommand{\dia}{$\diamond$}

\newcommand{\cM}{{\mathcal{M}}}
\newcommand{\cH}{\mathcal{H}}

\newcommand{\cP}{\mathcal{P}}

\newcommand{\bO}{\mathbf{O}}
\newcommand{\vol}{\mathrm{Vol}}

\newcommand{\tA}{\tilde{A}}
\newcommand{\bN}{\bar{N}}
\newcommand{\tZ}{\tilde{Z}}

\newcommand{\hA}{\hat{A}}

\begin{document}

\title{Why Polyhedra Matter in Non-Linear Equation Solving} 


\author{J.\ Maurice Rojas}
\thanks{ 
This research was partially supported by a grant {}from the 
Texas A\&M College of Science and NSF Grant DMS-0211458. } 
\address{
Department of Mathematics\\
Texas A\&M University\\
TAMU 3368\\
College Station, Texas 77843-3368\\
USA.}  
\email{rojas@math.tamu.edu \ {\tt http://www.math.tamu.edu/\~{}rojas} }  

\dedicatory{To my sister, Clarissa Amelia, on her 12$\thth$ birthday. }  

\begin{abstract} 
We give an elementary introduction to some recent polyhedral 
techniques for understanding and solving systems of multivariate 
polynomial equations. We provide numerous concrete examples and 
illustrations, and assume no background in algebraic geometry or convex 
geometry. Highlights include the following: 
\begin{enumerate} 
\item{A completely self-contained proof of an extension of 
Bernstein's Theorem. Our extension relates volumes of polytopes with the 
number of connected components of the complex zero set of a polynomial 
system, and allows any number of polynomials and/or variables.}  
\item{A near optimal complexity bound for computing {\bf mixed area} --- 
a quantity intimately related to counting complex roots in the plane} 
\end{enumerate} 

\end{abstract} 

\maketitle 

\section{Introduction} 
\label{sec:intro} 
In a perfect world, a scientist or engineer who wishes to solve a system of 
polynomial equations arising {}from some important application would simply 
pick up a book on algebraic geometry, look through the table of contents, and 
find a well-explained, provably fast algorithm which solves his or her 
problem. (Algebraic geometry began 2000 years ago as the study of polynomial 
equations, didn't it?) He or she would then surf the web to download a good 
(free) implementation which would run quickly enough to be useful. 

Once one stops laughing at how the real world compares, one realizes what is 
missing: the standard classical algebraic geometry texts (e.g., 
\cite{ega1,ega2,ega31,ega32,ega41,ega42,ega43,sga1,sga2,sga31,sga32,sga33,
sga41,sga42,sga43,sga45,sga5,sga6,sga7,mumford1,mumford2,hartshorne,
sha,gh}\footnote{In fairness, it should be noted that the major 
thrust of 20$\thth$ century algebraic geometry was understanding 
the {\bf topological} nature of zero sets of polynomials, rather than 
efficiently approximating the location of these zeros.})  
rarely contain algorithms and none contains a complexity analysis of any 
algorithm. Furthermore, one soon learns {}from experience that the specific 
structure 
underlying one's equations is rarely if ever exploited by a general 
purpose computational algebra package. 

Considering the ubiquity of polynomial equations in applications 
such as geometric modelling \cite{man98,goldman}, 
control theory \cite{sus98,nm99}, cryptography \cite{dod01}, 
radar imaging \cite{fh95}, learning theory \cite{vid97,vr02}, chemistry 
\cite{gh99,gat01}, game theory \cite{mcl97,roj97}, 
and kinematics \cite{can93,em99} (just to mention a few applications), 
it then becomes clear that we need an algorithmic theory of algebraic 
geometry that is practical as well as rigourous. One need only look at 
the active research in numerical linear algebra (e.g., eigenvalue problems 
for large sparse matrices) to see how far we are {}from a completely 
satisfactory theory for the numerical solution of general 
systems of multivariate polynomial equations. 

More recently, the introduction of algorithmic and combinatorial ideas 
has invigorated computational algebraic geometry. Here we give an elementary 
introduction to one recent aspect of computational algebraic geometry: 
polyhedral methods for solving systems 
of multivariate polynomial equations. The buzz-word for the 
cognicenti is {\bf toric varieties} \cite{tfulton,cox,frank}. However, rather 
than deriving algorithms 
{}from toric variety theory as an afterthought, we will begin directly with 
concrete examples and see how convex geometry naturally arises {}from solving 
equations. 

Simply put, polyhedral methods are a first step toward a new class of 
algorithms 
which adapt themselves to the intrinsic nature of the underlying system of 
equations. For the purposes of this paper, this means that our polynomial 
equations will be expressed as sums of monomial terms, and the techniques we 
describe will exploit the combinatorial structure of which monomial terms 
appear.

\begin{ex}
Suppose one has the following $3$ equations in $3$ unknowns $x$, $y$, and 
$z$:

\medskip
\noindent
\scalebox{.58}[1]{$c_{1,1}+c_{1,2}x+c_{1,3}y^2+c_{1,4}z^3+c_{1,5}x^5 y^6 z^7 + 
c_{1,6}x^6 y^7 x^5 + c_{1,7}x^7 y^5 z^6 + c_{1,8}x^8 y^9 z^9 + 
c_{1,9} x^{10} y^9 z^9 + c_{1,10} x^9 y^8z^9 +c_{1,11}x^9 y^{10} z^9+  
c_{1,12} x^9 y^9 z^{10}=0$}\\ 
\scalebox{.58}[1]{$c_{2,1}+c_{2,2}x+c_{2,3}y^2+c_{2,4}z^3+c_{2,5}x^5 y^6 z^7 + 
c_{2,6}x^6 y^7 x^5 + c_{2,7}x^7 y^5 z^6 + c_{2,8}x^8 y^9 z^9 + 
c_{2,9} x^{10} y^9 z^9 + c_{2,10} x^9 y^8z^9 +c_{2,11}x^9 y^{10} z^9+  
c_{2,12} x^9 y^9 z^{10}=0$}\\ 
\scalebox{.58}[1]{$c_{3,1}+c_{3,2}x+c_{3,3}y^2+c_{3,4}z^3+c_{3,5}x^5 y^6 z^7 + 
c_{3,6}x^6 y^7 x^5 + c_{3,7}x^7 y^5 z^6 + c_{3,8}x^8 y^9 z^9 + 
c_{3,9} x^{10} y^9 z^9 + c_{3,10} x^9 y^8z^9 +c_{3,11}x^9 y^{10} z^9+  
c_{3,12} x^9 y^9 z^{10}=0$,}

\medskip
\noindent
where the $c_{i,a}$ are any given complex numbers. One may reasonably 
guess that such a system of equations, being neither over-determined or under-determined, 
will have only finitely many roots $(x,y,z)\!\in\!\C^3$ with probability $1$, 
for any continuous probability distribution on the coefficient space 
$\C^{36}$. In fact, with probability $1$, the number of roots will 
always be the same number (cf.\ Theorems \ref{thm:disc} and \ref{thm:smoothkush} of Section 
\ref{sub:smooth}). What then is this ``generic'' number of roots? 

Noting that the maximum of the sum of the exponents in any summand of the 
first, second, or third equation is $28$ (i.e., our polynomials each have {\bf total degree} 
$28$), a classical theorem of B\'ezout \cite[Ex.\ 1, Pg.\ 198]{sha} gives us an upper bound 
of $21952\!=\!28^3$. A slightly more refined variant which uses 
degrees with respect to different variables, the {\bf multi-graded} 
version of B\'ezout's Theorem \cite{ms}, yields a sharper upper 
bound of $6000\!=\!6\cdot 10^3$.

However, the true generic number of roots is $\pmb{321}$. This number was 
calculated 
by using the correct concept in our setting: the convex hulls\footnote{ 
Recall that a set $B\!\subseteq\!\Rn$ is {\bf convex} iff for all 
$x,y\!\in\!B$, the line segment connecting $x$ and $y$ is also contained in 
$B$. The {\bf convex hull} of $B$, $\conv(B)$, is then simply the smallest 
convex set containing $B$, and the computational complexity of convex 
hulls of finite point sets is fairly well-understood \cite{preparata}. } 
of the exponent vectors (also known as the {\bf Newton polytopes}) of our 
polynomials.
In this case, all our Newton polytopes are identical, and the volume (suitably 
normalized) of any one 
serves as the correct generic number of complex roots. \\
\epsfig{file=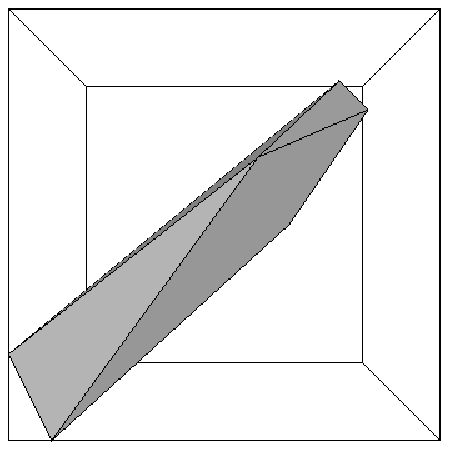,height=.68in} 
\epsfig{file=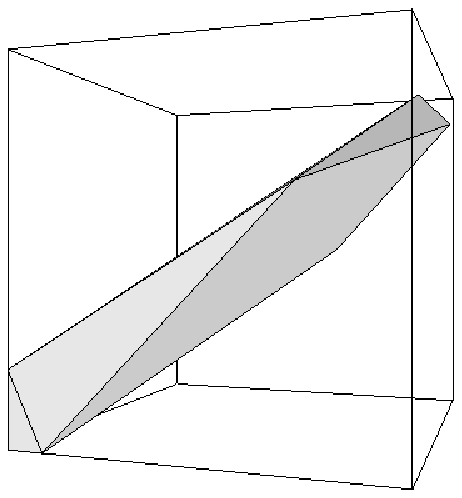,height=.68in} 
\epsfig{file=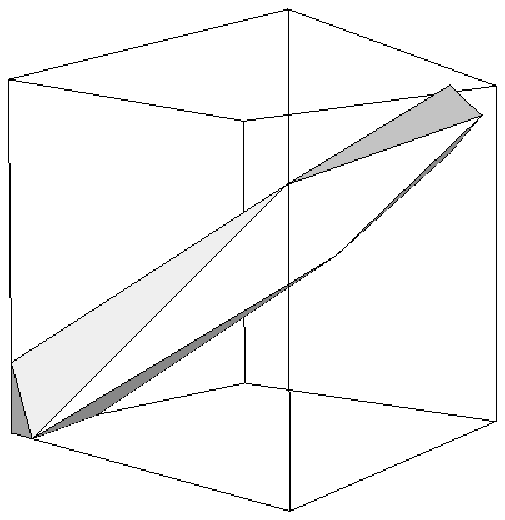,height=.68in} 
\epsfig{file=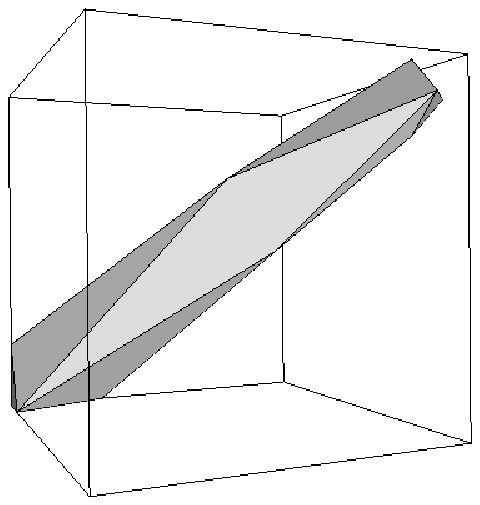,height=.68in} 
\epsfig{file=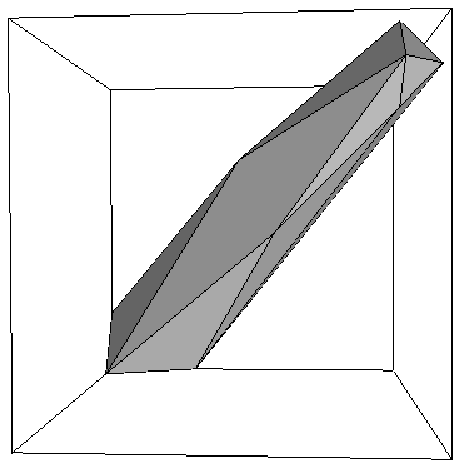,height=.68in} 
\epsfig{file=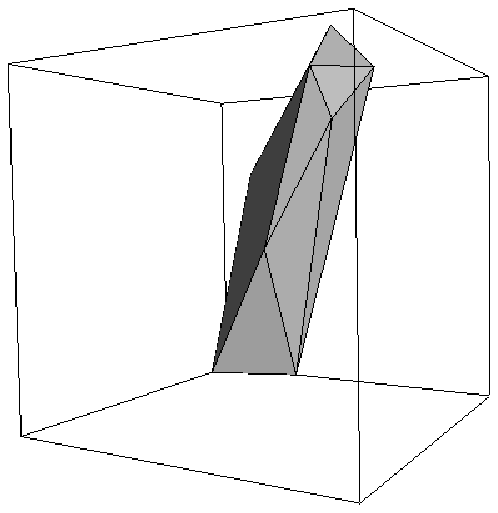,height=.68in} 
\epsfig{file=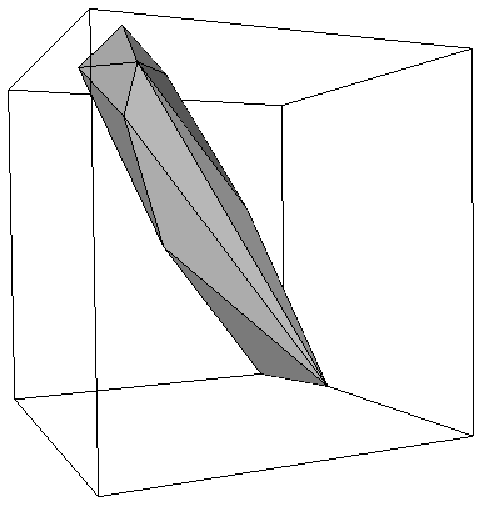,height=.68in}\\
\mbox{}\hfill{\footnotesize {\sc Figure 1.} {\em Several views of the Newton 
polytope shared by our three equations}}\hfill\mbox{} 
 
The key idea to keep in mind is that the complexity of 
solving a system of polynomial equations, or even approximating a 
{\bf single} root, depends strongly on the {\bf total} number of complex roots. Since one does 
not usually know this number a priori, the algorithm one uses ultimately makes implicit use 
of some upper bound on this number, usually one of the three we just 
mentioned. So, much as how our preceding comparison of bounds turned out, algorithms which 
take monomial term structure into account are preferable over those that do not. \dia 
\end{ex} 

A natural question, especially relevant to geometric modelling, then arises: 
Is there an analogous theory for systems of equations expressed in other 
bases? In particular, the systems of equations arising {}from $B$-splines are 
expressed in the so-called Bernstein-Bezier basis which uses sums of terms 
like $\prod_i(1-x_i)^{j_i}x^{k_i}_i$. The short answer is that an analogous theory for such 
bases does not yet exist. However, the philosophies of fewnomial theory 
\cite{few,tri,ari}, straight-line programs \cite{add,jkss}, not to mention polyhedral 
methods \cite{convex,hs,li97,gcp,verschelde,emirispan,mcdonald,high}, are bringing us closer 
to a theory that can efficiently handle such questions. 
 
We now outline the main results we explain in our paper. The 
second result below is new, while the first is older. 
(Related earlier results will be reviewed throughout this paper.) 
However, we emphasize that the statement and proof of the first 
result has been considerably simplified, we provide many more illustrations and 
examples than what is usually found in the earlier literature (e.g., 
\cite{bernie,kush1,kush2,tfulton,hs}), and we have made an effort to keep all 
prerequisites explicit and contained in this paper. 

\medskip

\noindent
\scalebox{.75}[1]{
\begin{minipage}[t]{3.125in} 
\begin{nota} 
Let $\bO$ denote the origin in 
\scalebox{.95}[1]{$\Rn$ and let $e_1,\ldots,e_n$ denote 
the standard basis vectors} 
\scalebox{.94}[1]{
\scalebox{1}[.8]{$\begin{bmatrix} 1 \\ 0 \\ \vdots \\ 0 \end{bmatrix},
\ldots,\begin{bmatrix} 0 \\ \vdots \\ 0 \\ 1 \end{bmatrix}$}$\in\!\Rn$. 
Also, for any $B\!\subseteq\!\Rn$, let $\conv(B)$} 
\scalebox{1.03}[1]{denote the smallest convex set containing $B$. Also,} 
\scalebox{1}[1]{we let $\vol(\cdot)$ denote the usual $n$-dimensional 
volume} 
\scalebox{.9}[1]{in $\Rn$, renormalized so that $\vol(\conv(\{\bO,
e_1,\ldots,e_n\})\!=\!1$.} Finally, we will abuse notation slightly by 
setting $\vol(A)\!:=\!\vol(\conv(A))$ whenever $A$ is a finite subset of 
$\Rn$. \dia 
\end{nota} 
\end{minipage}} 
\hspace{.1in}
\scalebox{.75}[1]{
\begin{minipage}[t]{3.125in} 
\begin{nota} 
For any $c\!\in\!\Cs\!:=\C\!\setminus\!\{0\}$ and
$a\!=\!(a_1,\ldots,a_n)\!\in\!\Zn$, let
$x^a\!:=\!x^{a_1}_1\cdots x^{a_n}_n$ and call $c x^a$ a
{\bf monomial term}. Also, for any polynomial of the form
$f(x)\!:=\!\sum_{a\in A} c_ax^a$, we call $\supp(f)\!:=\!\{a \; | \;
c_a\!\neq\!0\}$ the {\bf support} of $f$, and define
$\newt(f)\!:=\!\conv(\supp(f))$ to be the {\bf Newton polytope} of $f$.
We will assume henceforth that $F\!:=\!(f_1,\ldots,f_k)$ where, for all $i$,
$f_i\!\in\!\C[x_1,\ldots,x_n]$ and $\supp(f_i)\!\subseteq\!A_i$. 
We call $F$ a {\bf $\pmb{k\times n}$ polynomial system} (over $\C$) {\bf 
with support} $\pmb{(A_1,\ldots,A_k)}$. Finally, we let $Z_\C(F)$ denote 
\scalebox{.93}[1]{the set of 
$x\!\in\!\Cn$ with $f_1(x)\!=\cdots =\!f_k(x)\!=\!0$. \dia } 
\end{nota} 
\end{minipage}} 

\medskip
{{\sc Theorem 1.} (Special Case (full version in Sec.\ \ref{sec:direct})) 
\label{thm:complex} 
{\em 
Following the notation above, the number of connected 
components of $Z_\C(F)$ is no more than $\vol(B)$, where 
$B\!:=\!\{\bO,e_1,\ldots,e_n\}\cup \bigcup^k_{i=1} A_i$.
In particular, if the number of complex 
roots of $F$ is finite, then it is no more than $\vol(B)$. }} 

As might be expected, a sharper estimate on the generic number of 
complex roots comes at a price: the resulting formula is more 
difficult to evaluate. However, one can get an explicit and optimal 
complexity estimate for the case of a pair of bivariate equations. 

{{\sc Theorem 2.} 
\label{thm:complex2} 
{\em 
Following the notation of Theorem \ref{thm:complex}, 
suppose $k\!=\!n\!=\!2$. Then the generic number of complex roots of 
a polynomial system $F\!=\!(f_1,f_2)$ with support $(A_1,A_2)$ 
can be computed within $O(b\bN+\bN\log \bN)$ bit 
operations, and takes at least $\Omega(b\bN)$ bit operations in the worst case, 
where $\bN$ denotes the sum of the cardinalities of the $A_i$, and 
$b$ is the maximum bit-length of any coordinate of any $A_i$.} } 

Numerous examples of our two main theorems will appear as we 
review some of the background necessary for the applications  
and proofs of our theorems.  Theorem \ref{thm:complex} is proved three times: the 
simplified version above is proved in Section \ref{sec:poly} and then 
two proofs of the full version appear respectively in Sections \ref{sec:multi} 
and \ref{sec:mixed}. The last proof uses the main combinatorial 
construction detailed in this paper: the {\bf mixed subdivision} of 
$n$ polytopes in $n$ dimensions. Theorem 2 is proved in 
Section \ref{sec:mixed} as a simple consequence of mixed subdivisions 
in the plane. 

\section{From Binomial Systems to Volumes of Pyramids} 
\label{sec:binom} 
Perhaps the best and simplest place to begin to understand the 
connection between polytopes and polynomials is the special case of  
{\bf binomial} systems, i.e., polynomial systems where each polynomial 
has exactly $2$ monomial terms. For such systems, there is an 
immediate connection to linear algebra over the integers. 
\begin{ex} 
\label{ex:binom1}
Suppose we want to find all the complex solutions of
\begin{eqnarray*}
xy^7z^7w^4 & = & c_1\\
x^6y^4z^9w^6 & = & c_2\\
x^2y^3z^2w^6 &  = & c_3\\
x^6y^4z^8w^5 & = & c_4
\end{eqnarray*} 
where the $c_{i,j}$ are given nonzero complex numbers. Note in particular 
that this implies that any root of our system must satisfy $xyzw\!\neq\!0$. 
A particularly elegant trick we'll generalize shortly is the following: 
Consider the $4\times 4$ matrix 
$E\!:=${\small \scalebox{1}[.7]{$\begin{bmatrix}1 & 7 & 7 & 4\\ 
6 & 4 & 9 & 6 \\
2 & 3 & 2 & 6 \\
6 & 4 & 8 & 5
\end{bmatrix}$}} whose $i\thth$ row vector is the exponent vector of 
the $i\thth$ equation above. Then multiplying and dividing the equations 
above is easily seen to be equivalent to performing row operations on $E$. 
For example, doing a pivot operation to zero out all but the top entry of 
the first column of $E$ is just the computation of the matrix 
factorization {\small \scalebox{1}[.7]{$\begin{bmatrix} \mbox{} \ \ 1 & 0 & 0 & 0 \\
-6 & 1 & 0 & 0 \\
-2 & 0 & 1 & 0 \\
-6 & 0 & 0 & 1 \end{bmatrix} 
\begin{bmatrix}1 & 7 & 7 & 4\\ 
6 & 4 & 9 & 6 \\
2 & 3 & 2 & 6 \\
6 & 4 & 8 & 5
\end{bmatrix}
=\begin{bmatrix} 
1 & 7 & 7 & 4\\
0 & -38 & -33 & -18 \\
0 & -11 & -12 & -2 \\
0 & -38 & -34 & -19 
\end{bmatrix}$}}, which is in turn equivalent to observing that\\
\begin{tabular}{rcllll}
Equation 1 is...  &   $x$ & $y^7$ & $z^7$ & $w^4$ & $=c_{1}$\\
$\left.\text{(Equation 2)}\right/\!\text{(Equation 1)}^6 
\text{ is...}$ & & $y^{-38}$ & $z^{-33}$ & $w^{-18}$ & $=c^{-6}_{1}c_{2}$\\
$\left.\text{(Equation 3)}\right/\!\text{(Equation 1)}^2
\text{ is...}$ & & $y^{-11}$ & $z^{-12}$ & $w^{-2}$ & $=c^{-2}_{1}c_{3}$\\
$\left.\text{(Equation 4)}\right/\!\text{(Equation 1)}^6
\text{ is...}$ & & $y^{-38}$ & $z^{-34}$ & $w^{-19}$ & $=c^{-6}_{1}c_{4}$ 
\end{tabular}\\ 
Note also that this new binomial system has exactly the same roots 
as our original system. (This follows easily {}from the fact that 
our left-most matrix above is invertible, and the entries of the inverse are 
all {\bf integers}.) So we can solve the last $3$ equations for 
$(y,z,w)$ and then substitute into the first equation to solve for 
$x$ and be done. \dia 
\end{ex} 

Note, however, that if we wish to complete the solution of our 
example above, we must continue to use row operations on 
$E$ that are {\bf invertible over the integers}.\footnote{While one could 
simply use rational operations on $E$, and thus radicals on our equations, 
this quickly introduces some unpleasant ambiguities regarding 
choices of $d\thth$ roots. Hence the need for our integrality restriction.} 
This can be done by performing a simple variant of Gauss-Jordan\footnote{Recall that 
Gaussian elimination is what one does to reduce a matrix to upper triangular form 
(see Definition \ref{dfn:hermite} below), while Gauss-Jordan elimination is what one 
does to reduce a matrix to diagonal form \cite{strang}.} elimination  
where one uses {\bf no} divisions. In essence, one  
uses elementary {\bf integer} row operations to {\bf minimize the absolute value} of the 
entries in a given column, instead of reducing them to zero. 
\begin{ex} 
\label{ex:binom2} 
Let us consider the lower right $3\times 3$ block of our last matrix 
{}from our last example. The variant of Gauss-Jordan elimination we 
propose is easily depicted below.

\begin{picture}(200,20)(15,20)
\put(0,23){\epsfig{file=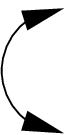,height=.16in}}
\put(5,20){{\small \scalebox{.7}[.8]{$\begin{bmatrix}-38 & -33 & -18 \\ 
-11 & -12 & -2 \\ 
-38 & -34 & -19 \end{bmatrix}$}}} 
\put(70,20){{\small \scalebox{.7}[.8]{$\longrightarrow$}}}
\put(87,30){\scalebox{.5}[.5]{$-3$}}
\put(92,11){\epsfig{file=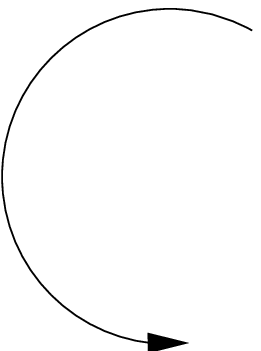,height=.32in}}
\put(98,22){\scalebox{.5}[.5]{$-3$}}
\put(105,22){\epsfig{file=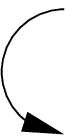,height=.16in}}
\put(110,20){{\small \scalebox{.7}[.8]{$ 
\begin{bmatrix}-11 & -12 & -18 \\ 
-38 & -33 & -18 \\ 
-38 & -34 & -19 \end{bmatrix}$}}} 
\put(175,20){{\small \scalebox{.7}[.8]{$\longrightarrow$}}}
\put(195,23){\epsfig{file=swap1.eps,height=.16in}}
\put(200,20){{\small \scalebox{.7}[.8]{$\begin{bmatrix}-11 & -12 & -18 \\ 
-5 & 3 & -12 \\ 
-5 & 2 & -13 \end{bmatrix}$}}} 
\put(265,20){{\small \scalebox{.5}[1.2]{$\longrightarrow etc\ldots 
\longrightarrow$}}}
\put(305,20){{\small \scalebox{.6}[.8]{$
\begin{bmatrix}-1 & -18 & 22 \\ 
\mbox{} \ 0 & 93 & -122 \\ 
\mbox{} \ 0 & -1 & -1 \end{bmatrix}$}}} 
\end{picture} 

\bigskip

Proceding with the same strategy on the lower $2\times 2$ block of the last 
matrix, and then continuing recursively, one finally obtains that our system of equations 
{}from Example \ref{ex:binom1} is equivalent to the following simpler system:
\begin{center}
\begin{tabular}{clllcllll}
$x$ &     &  & $w^{62}$  & $=$ & $c^{-3}_{1}$ & $c^{23}_2$    & $c^{11}_3$   & $c^{-26}_4$\\
    & $y$ &  & $w^{175}$ & $=$ & $c^{-8}_{1}$  & $c^{66}_{2}$ & $c^{31}_{3}$ & $c^{-75}_4$\\
    &     & $z$ & $w$ & $=$ & & $c^{1}_{2}$ & & $c^{-1}_{4}$ \\
 &  &  & $w^{215}$ & $=$ & $c^{-10}_{1}$ & $c^{82}_{2}$ & $c^{38}_3$ 
  & $c^{-93}_4$ 
\end{tabular}
\end{center} 
Note in particular that the exponent vectors above correspond exactly to 
the third and first matrices in the identity below: 
\begin{center}
{\small 
\scalebox{1}[.7]{$
\begin{bmatrix}-3 & 23 & 11 & -26\\ 
-8 & 66 & 31 & -75 \\
0 & 1 & 0 & -1 \\
-10 & 82 & 38 & -93 
\end{bmatrix}
\begin{bmatrix} 
1 & 7 & 7 & 4\\
6 & 4 & 9 & 6 \\
2 & 3 & 2 & 6 \\
6 & 4 & 8 & 5 
\end{bmatrix} 
=
\begin{bmatrix} 
1 & 0 & 0 & 62\\
0 & 1 & 0 & 175 \\
0 & 0 & 1 & 1 \\
0 & 0 & 0 & 215  
\end{bmatrix}$}}, 
\end{center} 
and that the left-most matrix $U$ has determinant $-1$ (and thus the entries 
of $U^{-1}$ are all integral). One should of course note that 
our system is now very easy to solve: since the roots of $w^d\!=\!c$ 
are $\left\{|c|^{1/d},|c|^{1/d}e^{2\pi i/d},\ldots,|c|^{1/d}e^{2\pi i(d-1)/d}\right\}$, 
we can simply substitute the $215$ resulting values for $w$ into the 
first, second, and third\\ 
\scalebox{.95}[1]{equations to solve for $x$, $y$, and $z$, 
thus finding all $215$ complex roots of our system. \dia} 
\end{ex} 

This motivates the following definition {}from 19$\thth$ century algebra. 
\begin{dfn} 
\label{dfn:hermite} 
Let $\Z^{n\times n}$ denote the set of all $n\times n$ matrices 
with all entries integral, and let $\gln(\Z)$ denote the 
set of all matrices in $\Z^{n\times n}$ with determinant $\pm 1$ 
(the set of {\bf unimodular} matrices). 
Recall that any $m\times n$ matrix $[u_{ij}]$ with 
$u_{ij}\!=\!0$ for all $i\!>\!j$ is called {\bf upper triangular}. 
Then, given any $M\!\in\!\Z^{n\times n}$, we call any 
identity of the form $UM=H$, with $H\!=\![h_{ij}]\!\in\!\Z^{n\times n}$ 
upper triangular and $U\!\in\!\gln(\Z)$, a {\bf Hermite factorization} 
of $M$. Also if, in addition, we have:  
\begin{enumerate}
\item{$h_{ij}\!\geq\!0$ for all $i,j$.}
\item{\scalebox{.96}[1]{for all $i$, if $j$ is the smallest $j'$ such that $h_{ij'}\!\neq\!0$ 
then $h_{ij}\!>\!h_{i'j}$ for all $i'\!\leq\!i$.}}  
\end{enumerate} 
then we call $H$ the {\bf Hermite normal form} of $M$. \dia 
\end{dfn} 
\begin{lemma} 
\label{lemma:unimod} 
\cite{unimod} 
For any $M\!\in\!\Z^{n\times n}$, a Hermite factorization can be 
computed within $O((n+h_M)^{6.376})$ bit operations, where 
$h_M\!:=\!\log(2n+\max\limits_{i,j}|m_{ij}|)$ and $M\!=\![m_{ij}]$. 
Furthermore, the Hermite normal form exists {\bf uniquely} for $M$, and can also 
be computed within the preceding bit complexity bound. \qed 
\end{lemma} 

By extending the tricks {}from our last examples, we can easily 
obtain the following lemma. 
\begin{lemma} 
\label{lemma:binom} 
Suppose $a_1,\ldots,a_n\!\in\Zn$ and $c_1,\ldots,c_n\!\in\!\Cs\!:=\!
\C\setminus\{0\}$. Let $E$ denote the $n\times n$ matrix whose 
$i\thth$ row is the vector $a_i$. Then the complex roots of the binomial 
system $F\!:=\!(x^{a_1}-c_1,\ldots,x^{a_n}-c_n)$ are exactly the 
complex solutions of the binomial system 
\begin{center}
\begin{tabular}{ccccccc}
$x^{u_{11}}_1$ & $\cdots$ & $x^{u_{1n}}_n$ & $=$ & $c^{u_{11}}_{1}$ & 
 $\cdots$ & $c^{u_{1n}}_{1}$ \\
  & $\ddots$ & $\vdots$ & $\vdots$ &  & $\vdots$ &  \\
 & & $x^{u_{nn}}_n$ & $=$ & $c^{u_{n1}}_{1}$ & 
 $\cdots$ & $c^{u_{nn}}_{1}$, \\
\end{tabular}
\end{center} 
where $[u_{ij}]E\!=\![h_{ij}]$ is any Hermite factorization of $E$. In particular, 
the complex roots of $F$ can be expressed explicitly as monomials in 
$\sqrt[h]{c_1},\ldots,\sqrt[h]{c_n}$, where $h\!:=\!\prod^n_{i=1}h_{ii}$. \qed 
\end{lemma} 

Letting $\Csn\!:=\!(\C\setminus\{0\})^n$, we then easily obtain the 
following corollary. 
\begin{dfn} 
Given any $k\times n$ polynomial system $F$, its {\bf Jacobian matrix} is 
$\jac(F)\!:=\!$\scalebox{1}[.7]{$\begin{bmatrix}
\frac{\partial f_1}{\partial x_1} & \cdots & \frac{\partial f_1}{\partial x_n}\\
\vdots & \ddots & \vdots \\
\frac{\partial f_k}{\partial x_1} & \cdots & \frac{\partial f_k}{\partial x_n}
\end{bmatrix}$}. We then say that a root $\zeta\!\in\!\Cn$  
of $F$ is {\bf degenerate} iff $\rank \; \jac(F)|_{x=\zeta}\!<\!k$, and 
{\bf smooth} otherwise. \dia 
\end{dfn}  
\begin{cor} 
\label{cor:binom} 
Suppose $F\!=\!(f_1,\dots,f_n)$ is any $n\times n$ binomial system and, for 
all $i$, $v_i$ is either vector defined by the difference of the exponent 
vectors of $f_i$. Then $F$ has only finitely many roots in $\Csn \Longrightarrow F$ 
has exactly $|\det M|$ many, where $M$ is the $n\times n$ matrix whose 
$i\thth$ row is $v_i$. In particular, the last quantity is exactly\\
$\vol(\{\bO,v_1,\ldots,v_n\})$. 

Also, every root of $F$ in $\Csn$ is non-degenerate iff 
$F$ has exactly $\vol(\{\bO,v_1,\ldots,v_n\})$ roots in $\Csn$. Finally, 
fixing the support of $F$, there is an algebraic hypersurface $\Delta$ in the 
coefficient space $\C^{2n}$ such that for all coefficient specializations 
specializations {\bf outside} of $\Delta$, $F$ has {\bf exactly} 
$\vol(\{\bO,v_1,\ldots,v_n\})$ roots in $\Csn$. \qed 
\end{cor} 
We illustrate the last portion of our corollary with the following 
example. 
\begin{ex} 
\label{ex:binom3} 
Let us find all $(c_{1,1},c_{1,2},c_{2,1},c_{2,2},c_{3,1},c_{3,2})\!\in\!
\C^6$ such that  
\begin{eqnarray*}
c_{1,1}x^2y^7z^5 & = & c_{1,2}\\
c_{2,1}x^4y^{14}z^{10} & = & c_{2,2}\\
c_{3,1}x^8y^{10}z^{14} &  = & c_{3,2}
\end{eqnarray*}
has infinitely many solutions in $(\Cs)^3$. In particular, by Lemma 
\ref{lemma:binom}, the roots of our system are exactly those of 
\begin{eqnarray*}
x^2\; y^7 \; z^5 & = & \frac{c_{1,2}}{c_{1,1}}\\
\mbox{} \; \ y^{18}z^{6} & = & \left(\frac{c_{1,2}}{c_{1,1}}\right)^{4}
\left(\frac{c_{3,2}}{c_{3,1}}\right)^{-1}\\
1 &  = & \left(\frac{c_{1,2}}{c_{1,1}}\right)^{-2}\frac{c_{2,2}}{c_{2,1}}
\end{eqnarray*}
Clearly then, our system has infinitely roots in $(\Cs)^3$ 
iff $c^{2}_{1,2}c_{2,1}\!=\!c_{2,2}c^2_{1,1}$ (and no roots whatsoever if 
$c^{2}_{1,2}c_{2,1}\!\neq\!c_{2,2}c^2_{1,1}$). So in this example, we can 
take 
\[\Delta\!=\!\left\{(c_{1,1},c_{1,2},c_{2,1},c_{2,2},c_{3,1},c_{3,2})\!\in\!
\C^6 \; | \; c^{2}_{1,2}c_{1,2}\!=\!c_{2,2}c^2_{1,1}\right\}. \text{ \dia} \]  
\end{ex} 

We conclude this section with a similar result for a slightly 
more complicated class of polynomial systems. 
\begin{dfn} 
\label{dfn:unmixed} 
Let $F$ be any $k\times n$ polynomial system with support $(A_1,\ldots,A_k)$. 
Then we say that $F$ {\bf is of type $\pmb{(m_1,\ldots,m_k)}$} iff 
$\#A_i\!=\!m_i$ for all $i$. Also, we that $F$ is {\bf unmixed} iff 
$A_1\!=\!\cdots\!=\!A_k$. Finally, writing $f_i(x)\!=\!\sum_{a\in A_i} 
c_{i,a}$ for all $i$, we say a property $\cP$ regarding $F$ holds 
{\bf generically} iff there is an algebraic hypersurface $\cH\!\subset\!
\C^{\sum_i \#A_i}$ such that $(c_{i,a} \; \; | \; \; i\!\in\!\{1,\ldots,n\} 
 \ , \ a\!\in\!A_i)\!\in\!\C^{\sum_i \#A_i}\setminus\cH \Longrightarrow 
\cP$ holds. \dia 
\end{dfn} 
\begin{prop} 
\label{prop:easy} 
Following the notation above, if $\cP_1,\ldots,\cP_\ell$ are properties of 
$F$ that hold generically, then their conjunction $\cP_1\wedge \cdots\wedge \cP_\ell$ 
holds generically as well. \qed 
\end{prop} 
\begin{cor} 
\label{cor:simplex}  
Given any unmixed $n\times n$ polynomial system $F\!=\!(f_1,\ldots,f_n)$ 
of type $(n+1,\ldots,n+1)$, let $A$ be the support of any $f_i$. 
Then $F$ either has exactly $\vol(A)$ roots in $\Csn$, no roots in 
$\Csn$, or infinitely many roots in $\Csn$. Furthermore, the first possibility 
holds generically and implies that all the roots of $F$ in $\Csn$ are non-degenerate. 
Finally, however many roots $F$ has in $\Csn$, they can always be expressed explicitly as 
monomials in $\vol(A)\thth$-roots of linear combinations of the coefficients of $F$ 
(and possibly some additional free parameters).  
\end{cor} 

\noindent
{\bf Proof:} By Gauss-Jordan elimination, $F$ is equivalent to a binomial 
system (i.e., one considers the monomials of $F$ as new variables, 
thus obtaining a linear system
that we can place into reduced 
row echelon form). So by Corollary \ref{cor:binom}, and 
some additional care with the Hermite normal when $F$ has infinitely many roots, 
we are done. \qed  

Corollary \ref{cor:simplex} will be the cornerstone of our proof of 
the special case of Theorem \ref{thm:complex} 
where $k\!=\!n$ and $F$ is unmixed (also known as {\bf Kushnirenko's Theorem}). Note in 
particular that any Newton 
polytope {}from a polynomial system as in Corollary \ref{cor:simplex}, 
when $\vol(A)\!>\!0$, is an $n$-simplex in $\Rn$, i.e., the $n$-dimensional 
analogue of a $3$-dimensional pyramid with a triangular base. 

\section{Subdividing Polyhedra and \mbox{Kushnirenko's} Theorem} 
\label{sec:poly} 
Here we prove the following central result which gives a strong 
connection between polytope volumes and the number of complex roots of 
polynomial systems. 
\begin{thm}[Kushnirenko's Theorem]   
\label{thm:kush} 
Following the notation of Definition \ref{dfn:unmixed} of 
Section \ref{sec:binom}, suppose $F\!=\!(f_1,\ldots,f_n)$ is 
$n\times n$ polynomial system with $\supp(f_i)\!\subseteq\!A$ for all $i$. 
Then $F$ has only finitely many roots in $\Csn \Longrightarrow F$ has 
$\leq\!\vol(A)$ roots in $\Csn$. Furthermore, for fixed $A$, 
$F$ generically has exactly $\vol(A)$ roots in $\Csn$. 
\end{thm} 
\noindent
This result is originally due to Anatoly Georievich Kushnirenko \cite{kush1, 
kush2}. So while we certainly claim no originality in our proof, we have 
strived to simplify the known proofs {}from the literature and include 
as much of the necessary background as reasonably possible.  

Before laying the technical foundations for our proof, let us 
first see a concrete illustration of the main ideas. 
In essence, one proves Kushnirenko's Theorem by deforming 
$F$ (preserving the number of roots along the way) into a 
{\bf collection} of simpler systems. Making this rigourous and 
efficient then provides a natural motivation for a new space (containing  
an embedded copy of $\Csn$) in which our roots will live. 
\begin{ex} 
\label{ex:2by2}
Consider the special case $n\!=\!2$ with 
\[ f_1(x,y)\!:=\!-2+x^2-3y+5x^7y^5+4x^6y^7 \] 
\[ f_2(x,y)\!:=\!3+2x^2+y+4x^7y^5+2x^6y^7. \] 
The Newton polygon boundary and support appear below:\\ 
\begin{picture}(133,100)(-70,-10)
\put(45,15){$\supp(f_1)\!=\!\supp(f_2)\!=\!A$}
\put(-10,0){\line(1,0){100}} \put(80,-7){$a_1$}
\put(0,-10){\line(0,1){100}} \put(-10,84){$a_2$}
\put(0,0){\circle*{3}} 
\put(20,0){\circle*{3}} 
 \put(20,0){\line(1,1){50}}
\put(70,50){\circle*{3}} 
 \put(70,50){\line(-1,2){10}}
\put(60,70){\circle*{3}} 
 \put(60,70){\line(-1,-1){60}}
\put(0,10){\circle*{3}} 
\end{picture}\\
According to Theorem \ref{thm:kush}, $F$ either has 
$\leq\!35$ roots in $(\Cs)^2$ or infinitely many. (The 
standard and multi-graded B\'ezout bounds respectively reduce to 
$169\!=\!13^2$ and $98\!=\!2\cdot 7\cdot 7$.) The true number of 
roots for our example turns out to be {\bf exactly} $35$, and these 
roots are all non-degenerate.  

To see why this is so, let us define a {\bf toric deformation} 
$\hat{F}\!:=\!(\hat{f}_1,\hat{f}_2)$ as follows: 
\[ \hat{f}_1(x,y,t)\!:=\!-2\pmb{t}+x^2-3y+5x^7y^5+4x^6y^7\pmb{t} \] 
\[ \hat{f}_2(x,y,t)\!:=\!3\pmb{t}+2x^2+y+4x^7y^5+2x^6y^7\pmb{t}, \] 
and $\hat{A}\!:=\!\supp(\hat{f}_1)\!=\!\supp(\hat{f}_2)\!=$
\scalebox{.8}[.8]{$\!\left\{\begin{bmatrix}
0 \\ 0 \\ 1 
\end{bmatrix},
\begin{bmatrix}
2 \\ 0 \\ 0 
\end{bmatrix},
\begin{bmatrix}
0 \\ 1 \\ 0 
\end{bmatrix},
\begin{bmatrix}
7 \\ 5 \\ 0 
\end{bmatrix},
\begin{bmatrix}
6 \\ 7 \\ 1 
\end{bmatrix}\right\}$}. (We will see later in Lemma \ref{lemma:lift} that 
most choices of powers of $t$ for each monomial would 
have worked just as well for proving Kushnirenko's Theorem for our example.) Note that 
$\hat{F}(1,x,y)\!=\!F(x,y)$.  Intuitively, one would expect $2$ equations in $3$ unknowns 
to generically define a curve (cf.\ Theorem \ref{thm:disc} and the 
Implicit Function Theorem from calculus), and this turns out to be the 
case for our example.  More to the point, the number of roots of $\hat{F}$ in 
$(\Cs)^2$ is {\bf constant} for all $t\!\in\!\C\!\setminus\!\Sigma$, where 
$\Sigma$ is a finite set not containing $1$.\footnote{ This crucial fact is elaborated a bit 
later in this section --- specifically, Lemma \ref{lemma:curve}.} 
So to show that $F$ has $35$ roots in $(\Cs)^2$, it suffices to 
show that the number of roots of $\hat{F}$ in $(\Cs)^2$ is $35$ for some 
suitable fixed $t$. At least initially, this seems no easier than counting 
the roots of $F$. 

The key trick then is to count something else which, for fixed $t$ 
sufficiently close to $0$, is easily provable to be the same as the number 
of roots of $\hat{F}$ in $(\Cs)^2$. This is where polyhedral subdivisions 
come into play almost magically. First, note that our new system is still 
unmixed but now has a $3$-dimensional Newton polytope:

\noindent 
\mbox{}\hspace{.5in}\epsfig{file=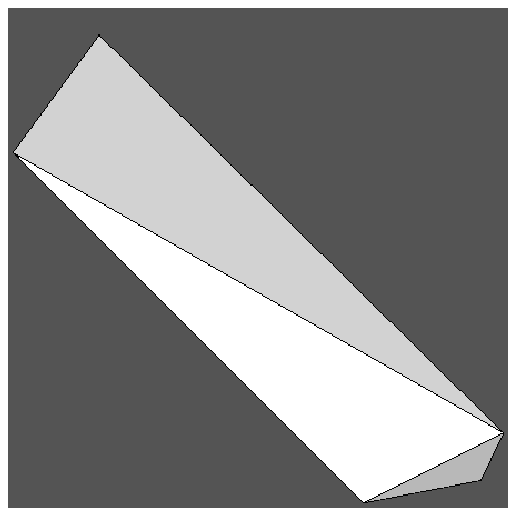,height=1.5in}
\hspace{1in}\epsfig{file=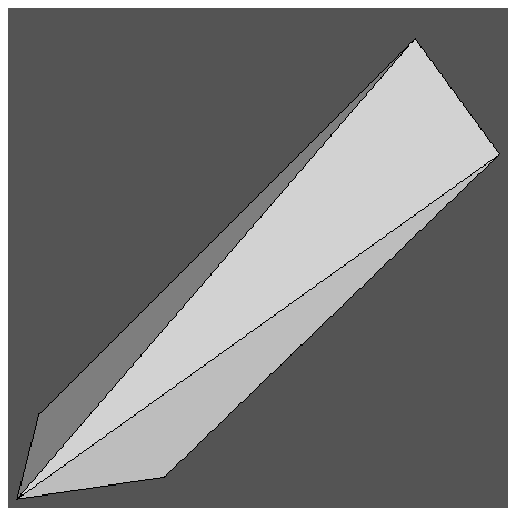,height=1.5in}\\
\mbox{}\hspace{.4in}{\small 
$\begin{matrix} 
\text{The {\bf Lower Hull} of $\conv(\hat{A})$}\\
\text{(View from }(0,0,t) \text{ with } t\!\ll\!0\text{)}
\end{matrix}$
\hspace{.75in}
$\begin{matrix} 
\text{The {\bf Upper Hull} of $\conv(\hat{A})$}\\
\text{(View from }(0,0,t) \text{ with }t\!\gg\!0\text{)}
\end{matrix}$}\\
Next, note that any root $(x,y,t)\!\in\!(\Cs)^3$ of $\hat{F}$ lies on a 
parametric curve of the form $C_{(x_0,y_0,w)}\!:=\!\{(s^{w_1}x_0,
s^{w_2}y_0,s^{w_3}) \; | \; s\!\in\!\C^*\}$. 
However, we will see momentarily that the set of $w\!\in\!\Zn$ for which 
the roots of $\hat{F}_t$ in $(\Cs)^3$ {\bf approach a  
$\pmb{C_{(x_0,y_0,w)}}$ as $\pmb{t\longrightarrow 0}$} is 
dictated by the face structure of $\conv(\hat{A})$, and 
{\bf all} the roots of $\hat{F}$ in $(\Cs)^3$ approach some 
$C_{(x_0,y_0,w)}$ as $t\longrightarrow 0$. \dia 
\end{ex} 

Let us pause now to review some basic convex geometric definitions. 
\begin{dfn}
A {\bf (closed) half-space (with (inner) normal 
$\pmb{a\!=\!(a_1,\ldots,a_n)}$)}, $H_a\!\subset\!\Rn$, is any set of the 
form $\{(y_1,\ldots,y_n)\!\in\!\Rn \; | \; a_1y_1+\cdots +a_ny_n\!\geq\!c\}$ 
for some real number $c$. A {\bf polyhedron} is any finite intersection 
of half-spaces. Also, a {\bf $\pmb{d}$-flat} in $\Rn$ is any translate of a 
$d$-dimensional subspace of $\Rn$. 
Finally, a {\bf (convex) polytope} in $\Rn$ is the convex 
hull of any finite point set in $\Rn$, and an {\bf $\pmb{n}$-simplex}  
is the convex hull of any $n+1$ points which do {\bf not} 
lie in an $(n-1)$-flat. \dia 
\end{dfn} 
\begin{dfn} 
For any $w\!:=\!(w_1,\ldots,w_n)\!\in\!\Rn$
and any compact set $B\!\subset\!\Rn$, we let $B^w$ ---
the {\bf face of $\pmb{B}$ with
(inner) normal $\pmb{w}$} --- be the set of all 
$y\!:=\!(y_1,\ldots,y_n)\!\in\!B$
minimizing the inner product $w\cdot y\!:=\!w_1y_1+\cdots+w_ny_n$. 
We call a face $Q$ of $B$ {\bf lower} (resp.\ upper) iff the last coordinate 
of any inner normal of $Q$ is positive (resp.\ lower).  
Also, the {\bf dimension of a face $Q$ of $B$}, $\dim Q$, is the 
dimension of the smallest flat containing $Q$. Finally, for any  
$d$-dimensional polytope, its faces of dimension $0$, $1$, $d-2$, 
$d-1$, or $d$ are respectively called {\bf vertices}, 
{\bf edges}, {\bf ridges}, {\bf facets}, or {\bf improper}. \dia 
\end{dfn} 

The most important thing we'll do with polytopes, after taking their 
faces, is to subdivide them. 
\begin{dfn} 
A {\bf subdivision} of a polytope $P$ is a collection of polytopes 
$\{Q_i\}^{\cN}_{i=1}$ called {\bf cells} satisfying the following conditions.  
\begin{enumerate} 
\item{$\bigcup\limits^n_{i=1} Q_i\!=\!P$ } 
\item{For any $w\!\in\!\Rn$ and $j\!\in\!\{1,\ldots,\cN\}$, we have 
$Q^w_j\!\in\!\{Q_i\}^{\cN}_{i=1}$ }
\end{enumerate} 
\scalebox{.95}[1]{In particular, if all the $Q_i$ are simplices then we say 
that $\{Q_i\}^{\cN}_{i=1}$ is a {\bf triangulation}.}\\ 
\scalebox{1}[1]{Finally, we call a 
cell $Q_j$ of $\{Q_i\}^{\cN}_{i=1}$ {\bf full-dimensional} iff 
$\dim Q_j\!=\!\dim P$. \dia} 
\end{dfn} 
\noindent
In particular, it is clear that one way to compute the volume of a 
polytope is to take any triangulation and add the volumes of all its 
full-dimensional cells. From basic linear algebra, we know that 
this reduces to a finite sum of absolute values of determinants 
of matrices of edge vectors. 

\addtocounter{ex}{-1}
\begin{ex}[Continued]
Let us now examine the lower hull of $\hat{A}$, projected 
onto the $(x,y)$-plane, and its inner lower facet normals.\\ 
\epsfig{file=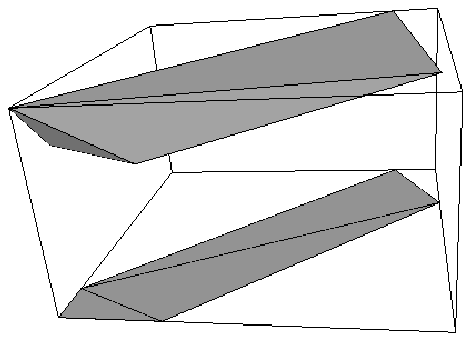,height=3in}
\scalebox{2}[2]{
\begin{picture}(133,80)(3,-20)
\put(1,2){\scalebox{.2}[.3]{$w\!=\!(1,2,2)$}}
\put(19,18){\scalebox{.3}[.3]{$w\!=\!(0,0,1)$}}
\put(-30,12){\begin{rotate}{-20}\epsfig{file=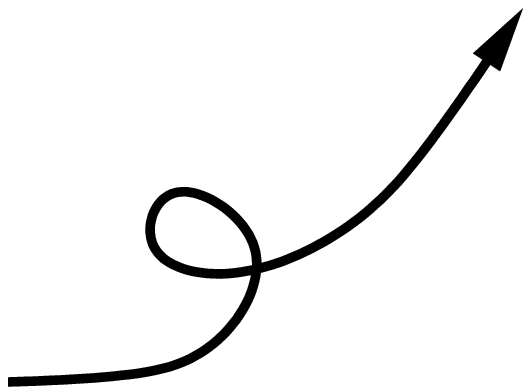,height=.3in}\end{rotate}} 
\put(-60,35){$\downarrow$} 
\put(37,44){\scalebox{.3}[.3]{$w\!=\!(4,-7,18)$}}
\put(0,0){\circle*{3}}
\put(0,0){\line(1,0){20}}
\put(0,0){\line(0,1){10}}
\put(20,0){\circle*{3}}
 \put(20,0){\line(1,1){50}}
\put(70,50){\circle*{3}}
 \put(70,50){\line(-1,2){10}}
\put(60,70){\circle*{3}}
 \put(60,70){\line(-1,-1){60}}
\put(0,10){\circle*{3}}
\put(20,0){\line(-2,1){20}}
\put(0,10){\line(7,4){70}}
\end{picture}}\\
In particular, the projections of the faces of the lower hull of $\hat{A}$ onto $A$ induce a
triangulation $\{Q_i\}$ of $\conv(A)$. 

Picking $w\!=\!(1,2,2)$ to examine the curves $C_{(x_0,y_0,w)}$, 
we see that $\hat{F}(s^{w_1}x_0,s^{w_2}y_0,s^{w_3})$ is 
exactly 
\[ s^2(-2+x^2_0-3y_0)+\text{Higher Order Terms in }s \]
\[ s^2(3+2x^2_0+y_0)+\text{Higher Order Terms in }s. \]
In particular, the $(x_0,y_0)\!\in\!(\Cs)^2$ which tend to a well-defined 
limit as $s\longrightarrow 0$ while satisfying $\hat{F}(s^1x_0,s^2y_0,s^2)\!=\!0$ 
must also satisfy 
\[ -2+x^2_0-3y_0 \] 
\[ 3+2x^2_0+y_0. \] 
(This follows easily from the Implicit Function Theorem upon observing that the roots of 
$(-2+x^2_0-3y_0, 3+2x^2_0+y_0)$ are all non-degenerate.) So by Corollary \ref{cor:simplex} of 
the last section, there are exactly $\vol(\{(0,0),(2,0),(0,1)\})\!=\!2$ 
such points. Put another way, the number of $(x_0,y_0)\!\in\!(\Cs)^2$  
for which $\hat{F}$ has roots in $(\Cs)^3$ approaching 
$C_{(x_0,y_0,(1,2,2))}$ as $t\longrightarrow 0$ is exactly $2$. 

Let us call the last system an {\bf initial term} 
system and observe that its Newton polytope is exactly the 
cell of $\{Q_i\}$ corresponding to $w\!=\!(1,2,2)$. 
Proceeding similarly with the other inner lower facet normals of 
$\hat{A}$, there are exactly $\vol(\{(0,1),(7,5),(6,7)\})\!=\!18$ curves 
of the form $C_{(x_0,y_0,(4,-7,18))}$, and exactly\\ 
$\vol(\{(2,0),(0,1),(7,5)\})\!=\!15$ curves of the form 
$C_{(x_0,y_0,(0,0,1))}$, approached by roots of $\hat{F}$ in $(\Cs)^3$  
as $t\longrightarrow 0$. Also, the last two initial term systems 
have Newton polytope respectively equal to the cell of $\{Q_i\}$ with 
inner lower facet normal $(4,-7,18)$ or $(0,0,1)$. 

To conclude, note that $w$ not a multiple of 
$(1,2,2)$, $(4,-7,18)$, or $(0,0,1) \Longrightarrow$ the resulting initial term 
system has a Newton polytope of dimension $\leq\!1$. Since 
$C_{(x_0,y_0,w)}\!=\!C_{(x_0,y_0,\alpha w)}$ for any $\alpha\!\in\!\Z$ and 
$w\!\in\!\Z^3$, another application of Corollary \ref{cor:simplex} then 
tells us that we have found {\bf all} $C_{(x_0,y_0,w)}$ (with 
$(x_0,y_0)\!\in\!(\Cs)^2$ and $w\!\in\!\Z^3$) that are 
approached by roots of $\hat{F}$ in $(\Cs)^3$ as $t\longrightarrow 0$. 
Since there are $35\!=\!\vol(A)$ such curves, and since they don't intersect 
at any fixed $t$, this implies that $\hat{F}_t$ has exactly $35$ roots in 
$(\Cs)^2$ for any $t$ with $|t|$ sufficiently small. So, assuming every root 
of $\hat{F}$ in $(\Cs)^3$ converges to some $C_{(x_0,y_0,w)}$ as $t\longrightarrow 0$, 
$F$ has exactly $35$ roots and we are done. \dia 
\end{ex} 

The preceding argument can be made completely general (not to mention 
rigourous) with just a little more work. In particular, we can prove our last 
assumption by constructing a space in which the roots of $\hat{F}$ all converge 
to well-defined limits as $t\longrightarrow 0$. This is one of the main 
motivations behind {\bf toric varieties}, which provide a useful and elegant way to 
compactify $\Csn$. 
\begin{rem} 
The idea of examining the behavior of the zero set of $F$ along a monomial 
curve is really not so far-fetched. In many practical computations (see, e.g., 
\cite{lw91,li97}), some roots of polynomial systems depending on parameters 
tend to have unbounded coordinates, even if the parameters are all finite. So 
monomial curves are a natural approach to formalize what it means for a root to 
``approach infinity.'' \dia 
\end{rem} 

\subsection{Enter Toric Varieties Corresponding to Point Sets}\mbox{}\\ 
\label{sub:toric} 
Before going further, let us first give a more succinct definition of 
initial term systems and formalize our constructions of $\hat{A}$ and 
$\hat{F}$.
\begin{dfn} 
For any $f\!\in\!\C[x_1,\ldots,x_n]$ of the form $\sum_{a\in A} c_ax^a$, 
let its {\bf initial term polynomial with respect to the 
weight $\pmb{w}$} be $\init_w(f)\!:=\!\sum_{a\in A^w} c_a x^a$. \dia 
\end{dfn} 
\begin{dfn} 
\label{dfn:lift} 
Given any $A\!\subset\!\Zn$, a {\bf lifting function} for $A$ is any  
function $\omega : A \longrightarrow \Rn$ and we let $\hat{A}\!:=\!\{(a,\omega(a)) 
\; | \; a\!\in\!A\}$. Also, letting 
$\pi : \R^{n+1} \longrightarrow \Rn$ denote the natural projection 
which forgets the last coordinate, we call 
$A_\omega\!:=\!\left\{\conv(\pi(\hat{A}^w)) \; | \; w\!\in\!\Rn\!\setminus\!\{\bO\}\right\}$ 
the {\bf subdivision 
of $\pmb{\conv(A)}$ induced by $\pmb{\omega}$}. Finally, we say that $\omega$ is a 
{\bf generic lifting} iff $A_\omega$ is a triangulation of $\conv(A)$. \dia
\end{dfn} 
\begin{dfn} 
Following the notation above, if 
we have in addition that $\omega(A)\!\subset\!\Zn$, then 
for any polynomial $f(x)\!=\!\sum_{a\in A} c_ax^a$, it's {\bf lift 
with respect to $\pmb{\omega}$} is the polynomial 
$\hat{f}(x,t)\!:=\!\sum_{a\in A} c_ax^at^{\omega(a)}$. 
Finally, the {\bf lift with respect to $\pmb{(\omega_1,\ldots,\omega_n)}$} of 
a $k\times n$ polynomial system $F\!:=\!(f_1,\ldots,f_k)$ is simply 
$\hat{F}\!:=\!(\hat{f}_1,\ldots,\hat{f}_k)$, where 
$\hat{f}_i$ is the lift of $f_i$ with respect to $\omega_i$ for all $i$. 
\dia  
\end{dfn} 
\begin{lemma}
\label{lemma:lift} 
Following the notation of Definition \ref{dfn:lift}, we have that for any 
fixed $A$, generic lifting functions occur generically. More precisely, 
there is a finite union, $\cH_A$, of proper flats in $\R^{\#A}$ such 
that $\omega(A)\!\in\!\R^{\#A}\setminus\cH \Longrightarrow \omega$ 
is a generic lifting for $A$. \qed 
\end{lemma} 
\noindent 

We will now refine the approach of Example \ref{ex:2by2} as follows: 
After building $\hat{A}$ and $\hat{F}$ via a generic lifting function, 
we will build a new point set $\tilde{A}$ and a space $Y_{\tilde{A}}$ with the following 
properties: 
\begin{enumerate} 
\item{$Y_{\tilde{A}}$ is compact} 
\item{There is an $h$-to-$1$ map from $\Csn$ to a dense open subset of 
$Y_{\tilde{A}}$, for some positive integer $h$.} 
\item{$\hat{F}$ has a well-defined complex zero set $\tZ$ in 
$Y_{\tilde{A}}$.}  
\item{There is a natural map $\pi : Y_{\tilde{A}} \longrightarrow 
\Pro^1_\C$, where $\Pro^1_\C\!=\!\C\cup\{\infty\}$ is the usual 
complex projective line, such that  
for all $t_0\!\in\!\Cs$, $h\#(\pi^{-1}(t_0)\cap \tZ)$ is 
exactly the number of roots of $\hat{F}$ in $\Csn$ with 
$t$-coordinate $t_0$.}    
\end{enumerate}
Our proof of Kushnirenko's Theorem will then focus instead 
on (a) showing that $\#(\pi^{-1}(1)\cap\tZ)\!=\!\#(\pi^{-1}(0)\cap\tZ)$ generically, and 
(b) computing $\#\pi^{-1}(t)$ {\bf at} $t\!=\!0$ to avoid the use of limits. We've actually 
already seen an example of (b), from an elementary point of view, in Example 
\ref{ex:2by2} of the last section. So let us now elaborate the framework needed for (a).  

\begin{dfn} 
\label{dfn:toric} 
Let $N\!:=\!\#A$ and let \[\Pro^N_\C\!:=\!\{[x_0 : \cdots : x_N] \; | \; 
x_0,\ldots,x_N\!\in\!\C \text{ not {\bf all} } 0\}\] denote 
complex $N$-dimensional projective space. 
Then, given any finite subset $A\!\subset\!\Zn$, we let 
$\varphi_A : \Csn \longrightarrow \Pro^{N-1}_\C$ --- the 
{\bf generalized Veronese map} --- be the map defined by 
$x \mapsto [x^{a_1}:\cdots:x^{a_N}]$, where $\{a_1,\ldots,a_N\}\!=\!A$.  
We then let $Y_A$ --- the {\bf toric variety corresponding to the point set 
$\pmb{A}$} --- denote the closure of $\varphi_A\left(\Csn\right)$ in $\Pro^{N-1}_\C$. 
\dia
\end{dfn} 

Being a closed subset of a compact space, we thus see that $Y_A$ is 
compact as a topological space and this will be important later for 
guaranteeing that certain limits of curves exist. However, one may wonder if $Y_A$ actually 
compactifies $\Csn$ in any reasonable way and what the closure above really means. Here's one 
way to make this precise. 
\begin{lemma} 
\label{lemma:veronese}
Following the notation of Definition \ref{dfn:toric}, 
let $[a_{i1},\ldots,a_{in}]\!=\!a_i$ for all $i$. 
Also let $E$ (resp.\ $\bar{E}$) be the $N\times n$ (resp.\ 
$N\times (n+1)$) matrix whose $i\thth$ row is $(a_{i1},\ldots,a_{in})$ 
(resp.\ $(a_{i1},\ldots,a_{in},1)$). Finally, let 
$H$ be the Hermite normal form of $E$, let $\bar{U}\bar{E}\!=\!\bar{H}$ be any Hermite 
factorization of $\bar{E}$, and let $\bar{u}_i$ (resp.\ $h$) denote the 
$i\thth$ row of $\bar{U}$ (resp.\ the product of the 
diagonal elements of $H$). Then 
$ Y_A\!=\!\left\{[p_1:\cdots : p_N]\!\in\!\Pro^{N-1}_\C \; | \;  
p^{\bar{u}^+_{r+1}}\!=\!p^{\bar{u}^-_{r+1}},\ldots,p^{\bar{u}^+_N}\!=\!
p^{\bar{u}^-_N}\right\}$, 
where $r$ is the rank of $\bar{H}$ and, for all $i$, $\bar{u}^+_i-\bar{u}^-_i\!=\!\bar{u}_i$  
and $\bar{u}^{\pm}_i\!\in\!\Rn_+$. 
Furthermore, $\varphi_A$ is an $h$-to-$1$ map, i.e., $\#\varphi^{-1}_A(p)\!=\!h$ for all 
$p\!\in\!\varphi_A\left(\Csn\right)$. \qed 
\end{lemma} 
\noindent 
The $p_i$ above are sometimes called {\bf toric coordinates}. The proof of 
Lemma \ref{lemma:veronese} is a routine application of the Hermite normal form we 
introduced in the last section. Let us see an example of $Y_A$ now.
\begin{ex} 
Taking $A$ as in our last example, we obtain 
\[\varphi_A(x,y)\!=\![1:x^2:y:x^7y^5:x^6y^7] \ , \  
E\!=\!\text{\scalebox{.6}[.5]{$\begin{bmatrix} 0 & 0 \\ 2 & 0 \\ 0 & 1 \\ 7 & 5 \\ 6 & 7 
\end{bmatrix}$}} \ , \ 
\text{and } \bar{E}\!=\!\text{\scalebox{.6}[.5]{$\begin{bmatrix} 0 & 0 & 1\\ 2 & 0 & 1\\ 
0 & 1 & 1\\ 7 & 5 & 1\\ 6 & 7 & 1 \end{bmatrix}$}}.\] Using the {\tt ihermite} command 
in {\tt Maple}, we then easily 
obtain that $H\!=\!$\scalebox{.6}[.5]{$\begin{bmatrix} 1 & 0 \\ 0 & 1 \\ 0 & 0 \\ 
0 & 0 \\ 0 & 0 \end{bmatrix}$} is the Hermite normal form for $E$ and 
\scalebox{.6}[.5]{$\begin{bmatrix} 7 & -3 & -5 & 1 & 0\\ -1 & 0 & 1 & 0 & 0 \\ 
1 & 0 & 0 & 0 & 0\\ 15 & -7 & -10 & 2 & 0 \\ 9 & -3 & -7 & 0 & 1 
\end{bmatrix}$}$\bar{E}\!=\!$\scalebox{.6}[.5]{$\begin{bmatrix} 1 & 0 & 0 \\ 
0 & 1 & 0 \\ 0 & 0 & 1\\ 0 & 0 & 0 \\ 0 & 0 & 0 
\end{bmatrix}$} is a Hermite factorization for $\bar{E}$.  
So Lemma \ref{lemma:veronese} tells us that our $Y_A$ here can also be defined as the zero 
set in $\Pro^4_\C$ of the following collection of binomials: 
\[ \left\langle \ p^{15}_1p^2_4-p^7_2p^{10}_3 \ , \ p^9_1p_5-p^3_2p^7_3 \ \right\rangle. \] 
(Note, for instance, that 
$1^{15}(x^7y^5)^2-(x^2)^7y^{10}\!=\!x^{14}y^{10}-x^{14}y^{10}\!=\!0$.) 
Furthermore, since $h\!=\!1\cdot 1\cdot 1\!=\!1$, our map $\varphi_A$ here is thus a bijection. \dia 
\end{ex} 

The most relevant combinatorial/geometric properties of $Y_A$ can 
be summarized as follows (see the companion tutorials \cite{cox,frank} in this volume 
for other aspects and points of view). 
\begin{dfn} 
\label{dfn:equiv} 
Following the notation above, let $N\!=\!\#A$ and $\{a_1,\ldots,a_N\}\!=\!A$. 
Then, for any face $Q$ of $\conv(A)$, the {\bf orbit} $O_Q$ is the subset 
\[\{[p_1:\cdots :p_N]\!\in\!Y_A \; | \; p_i\!=\!0 \Longrightarrow 
a_i\!\not\in\!Q\}.\] 
Also, for any $p\!\in\!O_Q$ with $Q$ a {\bf proper} face, 
we say that $p$ {\bf lies at toric infinity}. Finally, given any 
$f_1,\ldots,f_k\!\in\!\C[x_1,\ldots,x_n]$ of 
the form $f_i(x)\!=\!\sum_{a\in A} c_{i,a} x^a$ for all $i$, the {\bf zero set of $
\pmb{F\!=\!(f_1,\ldots,f_k)}$ in $\pmb{Y_A}$} is simply the set of all 
$[p_1:\cdots :p_N]\!\in\!Y_A$ with $\sum^N_{j=1}c_{i,a_j} p_j\!=\!0$ for all $i$. \dia 
\end{dfn} 
\noindent
Just as a polytope can be expressed as a disjoint union of the relative interiors 
of its faces, $Y_A$ can always be expressed as disjoint union of the 
$O_Q$. The lemma below follows routinely from Lemma \ref{lemma:veronese} and 
Definition \ref{dfn:equiv}. Recall that a {\bf (complex) algebraic set} is simply a 
subset of $\C^N$ or $\Pro^N_\C$ defined by the zero set of a system of 
polynomial equations. 
\begin{lemma} 
\label{lemma:face}
Following the notation of Definition \ref{dfn:equiv}, $O_Q$ is a dense open 
subset of a $d$-dimensional algebraic subset of $\Pro^{N-1}_\C$, where 
$d\!=\!\dim Q$. Also, $Y_A$ is the disjoint union 
\[\bigsqcup\limits_{Q \text{ a face of } \conv(A)} O_Q\] 
and \[Y_A\setminus\varphi_A\left(\Csn\right)\!=\!\bigsqcup
\limits_{Q \text{ a {\bf proper} face of } \conv(A)} O_Q.\] 
Finally, $F$ has a root in $O_{\conv(A)^w}$ iff $\init_w(F)$ has a root in 
$\Csn$. \qed 
\end{lemma} 
\noindent
Since all faces of $\conv(A)$ have a well-defined inner normal, 
Lemma \ref{lemma:face} thus gives a complete characterization of 
when a root of $F$ lies at toric infinity, as well as which piece 
of toric infinity. This is what will allow 
us to replace the cumbersome curves $C_{(x_0,y_0,w)}$ mentioned earlier in 
Example \ref{ex:2by2} with a single algebraic curve in $Y_A$. 

\subsection{The Smooth Case of Kushnirenko's Theorem}\mbox{}\\ 
\label{sub:smooth} 
Let us now review the last tool we'll need to start our proof of Kushnirenko's Theorem:  
Simplified characterizations of the {\bf $\pmb{A}$-discriminant} and {\bf Cayley trick}, 
and some basic facts on algebraic curves. 
\begin{dfn} 
\label{dfn:tjac} 
Given any $k\times n$ polynomial system $F$, 
the {\bf toric} Jacobian matrix of $F$ 
is $\toricjac(F)\!=\!$\scalebox{1}[.7]{$\begin{bmatrix}
x_1    & 0       & \cdots  & 0      \\
0      & x_2  &         &        \\
\vdots &         &         & \vdots \\  
       &         & \ddots  & 0      \\
0      & \cdots  & 0       & x_n 
\end{bmatrix}$}$\jac(F)$. Assuming $F$ is unmixed, we then say that $F$ 
{\bf has a degenerate root at $\pmb{p\!=\![p_1:\cdots:\!p_N]\!\in\!Y_A}$} 
iff $p$ is a root of $F$ and $\rank \; \toricjac(F)|_p\!<\!k$. 
We then let the {\bf discriminant variety}, $\Delta_{(\underset{k}{\underbrace{A,\ldots,A}})}$, 
denote the set of all 
\vspace{-.7cm}
\[(c_{1,a_1},\ldots,c_{1,a_N})\times \cdots \times 
(c_{k,a_1},\ldots,c_{k,a_N})\!\in\!(\C^N)^k\] such that $F$ has a degenerate root 
in $Y_A$. Finally, given any $k$-tuple of point sets from $\Rn$, $(A_1,\ldots,A_k)$, its 
{\bf Cayley configuration} is the point set \[\cay(A_1,\ldots,A_k)\!:=\!
(A_1\times (\underset{k-1}{\underbrace{0,\ldots,0}}))\cup 
(A_2\times e_{n+1}) \cup \cdots \cup (A_k\times e_{n+k-1})\!\subset\!\R^{n+k-1}. 
\text{ \dia } \]   
\end{dfn}  
\begin{ex} 
Returning to Example \ref{ex:2by2} one last time, consider the 
root $q\!=\![1:\sqrt{-1}:-1:0:0]\!\in\!Y_{\hat{A}}$ of 
\[ f_1(x,y)\!:=\!-2t+x^2-3y+5x^7y^5+4x^6y^7t \] 
\[ f_2(x,y)\!:=\!3t+2x^2+y+4x^7y^5+2x^6y^7t, \] 
where $\hat{A}\!:=\!\supp(\hat{f}_1)\!=\!\supp(\hat{f}_2)\!=$\scalebox{.8}[.8]
{$\!\left\{\begin{bmatrix}
0 \\ 0 \\ 1 
\end{bmatrix},
\begin{bmatrix}
2 \\ 0 \\ 0 
\end{bmatrix},
\begin{bmatrix}
0 \\ 1 \\ 0 
\end{bmatrix},
\begin{bmatrix}
7 \\ 5 \\ 0 
\end{bmatrix},
\begin{bmatrix}
6 \\ 7 \\ 1 
\end{bmatrix}\right\}$}. Note in particular that $q\!\in\!O_{(1,2,2)}$ and thus lies at the 
portion of toric infinity corresponding to the smallest triangular cell of 
$A_\omega$. The toric Jacobian matrix, in toric coordinates, is then 
\scalebox{.8}[1]{$\begin{bmatrix} 
2p_2 +35p_4+24p_5 & -3p_3 +25p_4+28p_5 & -2p_1 + 4p_5 \\
4p_2 +28p_4+12p_5 & p_3 + 20p_4+14p_5 & 3p_1+2p_5
\end{bmatrix}$}. Evaluating at $q$, our matrix then becomes 
\scalebox{.8}[1]{$\begin{bmatrix} 
2\sqrt{-1} & 3  & -2 \\
4\sqrt{-1} & -1 & 3
\end{bmatrix}$}, which clearly has rank $2$, so $q$ is a non-degenerate root. \dia 
\end{ex} 
\begin{thm} 
\label{thm:disc} 
Following the notation of Definition \ref{dfn:tjac}, there is a 
homogeneous polynomial $D_A\!\in\!\C[c_{a_1},\ldots,c_{a_N}]$ such that 
\[D_A(c_{a_1},\ldots,c_{a_N})\!\neq\!0 \Longrightarrow (c_{a_1},\ldots,c_{a_N})\!\not\in\!
\Delta_A,\] i.e., $\Delta_A$ is always contained in an algebraic hypersurface in 
$\C^N$. In particular, \[D_{\cay(\underset{k}{\underbrace{A,\ldots,A}})}\!\neq\!0 
\Longrightarrow (c_{1,a_1},\ldots,c_{1,a_N})\times \cdots \times
(c_{k,a_1},\ldots,c_{k,a_N})\!\not\in\!\Delta_{A^k}. \text{ \qed}\] 
\end{thm}
\begin{cor} 
\label{cor:0vol}
Suppose $F$ is an $n\!\times n$ polynomial system with support $(A_1,\ldots,A_n)$ and 
that there is an $(n-1)$-flat containing translates of $A_1,\ldots,A_n$. Then for fixed 
$(A_1,\ldots,A_n)$, $F$ generically has no roots in $\Csn$. In particular, in the 
unmixed case, $F$ generically has no roots in $Y_A$. \qed 
\end{cor}  
\noindent
Obtaining the discriminant of a system of equations via the discriminant of a 
single larger equation via the Cayley configuration is sometimes called the 
{\bf Cayley trick} \cite{gkz}. Theorem \ref{thm:disc} can actually be derived 
directly from our framework here via the {\bf toric resultant} (see, e.g.,  
\cite{mourrain} in this volume or \cite{emirispan}). However, for the sake of 
brevity we omit the proof.  
The final additional fact we'll need follows easily from the Implicit Function Theorem.
\begin{dfn}
If $X\subseteq\!\Pro^{N_1}_\C$ and 
$Y\!subseteq\!\Pro^{N_2}_\C$ are algebraic curves, then a {\bf morphism} 
$\psi : X \longrightarrow Y$ is simply a well-defined map of the form  
$[p_1:\cdots :p_{N_2}]\mapsto [\phi_1(p_1,\ldots,p_{N_2}):\cdots:
\phi_{N_1}(p_1,\ldots,P_{N_2})]$, where $\phi_1,\ldots,\phi_{N_2}$ are 
homogeneous polynomials of the same degree. \dia 
\end{dfn} 
\begin{lemma} 
\label{lemma:curve}
Suppose $C\!\subset\!\Pro^N_\C$ is a smooth algebraic curve (not necessarily 
connected) and $\psi : C \longrightarrow \Pro^1_\C$ is any morphism.  Then either 
$\#\psi(X)\!<\!\infty$ or $\psi(X)\!=\!\Pro^1_\C$. In the latter 
case, there is a positive integer $m$ and a finite set $\crit_\psi\!\subset\!\Pro^1_\C$, the 
{\bf critical values} of $\psi$, such that $\#psi^{-1}(t)\!=\!m \Longleftrightarrow  
t\!\in\!\Pro^1_\C\!\setminus\!\Sigma$. 

Finally, in the special case where $C$ is the zero set in $Y_{\tilde{A}}$ of an 
$n\times (n+1)$ polynomial system $\hat{F}(x_1,\ldots,x_n,t)$ with 
$\supp(\hat{f})\!\subseteq\!\hat{A}$ 
for all $i$, $\hat{A}\!\subseteq\!\tilde{A}$, 
and $\psi(\varphi_{\hat{A}}(x_1,\ldots,x_n,t))\!=\![1:t]$ for all $t\!\in\!\Cs$, 
we have that $t_0\!\in\!\C$ lies in $\crit_\psi \Longleftrightarrow (\hat{F},t-t_0)$ has a 
degenerate root in $Y_{\tilde{A}}$.   
\end{lemma}  
\noindent

\begin{thm}
\label{thm:smoothkush}
Following the notation of Theorem \ref{thm:kush}, fix $A$. 
Then $F$ generically has exactly $\vol(A)$ roots in $\Csn$, all 
of which are non-degenerate. 
\end{thm}

\noindent 
{\bf Proof of Theorem \ref{thm:smoothkush}:} Let $N\!:=\#A$ and $\{a_1,\ldots,a_N\}\!:=\!A$ as 
before and pick any generic lifting function $\omega$ with {\bf integral} range. 
Following the notation of Definition \ref{dfn:lift}, 
let us then define $\tilde{A}\!:=\!\hat{A}\cup\{(a,\omega(a)+1) \; | \; a\!\in\!A\}$ 
and order the coordinates $p\!=\![p_1:\cdots :p_{2N}]$ of $\Pro^{2N}_\C$ so 
that \[\varphi_{\tilde{A}}(x,t)\!=\![x^{a_1}t^{\omega(a_1)}:\cdots :x^{a_N}t^{\omega(a_N)}:
x^{a_1}t^{\omega(a_1)+1}:\cdots :x^{a_N}t^{\omega(a_N)+1}].\] Note in particular 
that $\hat{F}$ has a well-defined zero set in $Y_{\tilde{A}}$, as well as $Y_{\hat{A}}$: 
since the coordinates $p_1,\ldots,p_N$ of $Y_{\hat{A}}$ can be identified with an obvious 
subset of the coordinates of $Y_{\tilde{A}}$. 

We can now at last define our promised map $\pi : Y_{\tilde{A}} \longrightarrow 
\Pro^1_\C$ by $p\mapsto [p_1:p_{N+1}]$. Defining $\tZ$ (resp.\ $Z$) to be the 
zero set of $\hat{F}$ in $Y_{\tilde{A}}$ (resp.\ $F$ in $Y_A$), note that 
there is an isomorphism (an algebraic bijection) between 
$\pi^{-1}(1)\cap \tZ$ and $Z$ defined by $[p_1:\cdots :p_{2N}] \longleftrightarrow 
[p_1:\cdots:p_N]$. (This is easily checked since Lemma \ref{lemma:veronese} tells us 
that the binomials that define $Y_{\tA}$ are exactly those defining $Y_A$ union 
$\langle p_{N+1}p_2-p_1p_{N+2},\cdots,p_{N+1}p_N-p_1p_{2N}\rangle$.) 

Now note that $\pi$ also induces a natural morphism from $\tZ$ to $\Pro^1_\C$. Let $H$ be the 
Hermite normal form of $A$ and $h$ the product of the diagonal elements of $H$. Since the 
first $n$ columns of the Hermite normal forms of $A$ and $\tA$ are the 
same, Lemma \ref{lemma:veronese} then tells us that the number of roots of $F$ is exactly 
$h\#(\pi^{-1}(1)\cap \tZ)$. By applying Theorem \ref{thm:disc} to 
$(A,\ldots,A)$ (and Proposition \ref{prop:easy}) 
it thus suffices to show that $\pmb{h\#(\pi^{-1}(1)\cap \tZ)\!=\!\vol(A)}$ 
generically. 

Next, note that by construction, all the initial term systems of $F$ will be 
unmixed and have Newton polytopes of volume $0$. In particular, by Corollary 
\ref{cor:0vol}, any particular initial term system will generically have no roots. 
Similarly, by Corollary \ref{cor:simplex}, the initial term systems of $\hat{F}$ will have 
each have smooth zero set generically. 
So by Lemma \ref{lemma:face} and Proposition \ref{prop:easy}, it will be 
generically true that $F$ will have {\bf no} roots at toric infinity in 
$Y_A$, and all the roots of $\hat{F}$ at toric infinity in 
$Y_{\tilde{A}}$ will be non-degenerate. Furthermore, by applying Theorem \ref{thm:disc} to 
$(\hat{A},\ldots,\hat{A})$, we know that $\tZ$ is generically smooth. So by 
Proposition \ref{prop:easy} again, {\bf it thus suffices to show that 
[$\pmb{\tZ}$, $\pmb{Z}$, and $\pmb{\tZ\cap \text{(Toric Infinity in }Y_{\tilde{A}}\text{)}}$ 
smooth and $\pmb{Z\!\subset\!\varphi_{\tilde{A}}(\Csn)}$]   
$\pmb{\Longrightarrow h\#(\pi^{-1}(1)\cap \tZ)\!=\!\vol(A)}$.} 

So let us now assume the hypothesis of the last implication. 
By Lemma \ref{lemma:curve}, $Z$ (resp.\ $\tZ\cap
 \pi^{-1}(0)$) smooth $\Longrightarrow 1$ (resp.\ $0$) is not a critical value of 
$\pi$. Also, by the Implicit Function Theorem, the smoothness of $\tZ$ implies that 
$\pi(\hat{Z})$ contains a small open ball about $1$. So by the first part of 
Lemma \ref{lemma:curve}, $\pi(\tZ)\!=\!\Pro^1_\C$.

Clearly, $\Pro^1_\C$ remains path-connected even after a finite set of 
points is removed, so let $L$ be any continuous path connecting $0$ and $1$ in 
$\Pro^1_\C\!\setminus\!\crit(\pi|_{\tZ})$.  By the Implicit Function Theorem once more, and 
the fact that $L$ is compact (by virture of the compactness of $\Pro^1_\C$), we must have that 
$\#(\pi^{-1}(t)\cap\tZ)$ is constant on $L$. {\bf So we now need only show that 
$\pmb{h\#(\pi^{-1}(0)\cap\tZ)\!=\!\vol(A)}$.}

To conclude, note that $\tA$ and $\hA$ have the same lower hull, so Lemmata \ref{lemma:face} 
and \ref{lemma:veronese} then imply that $\pi^{-1}(0)\cap\tZ$ is nothing more than 
\[\left\{[p_1:\cdots:p_{2N}]\!\in\!Y_{\tilde{A}} \; \; \left| \; \;  
\sum\limits_{a_j\in Q} c_{i,{a_j}}p_j\!=\!0 \text{ for all } i\!\in\!\{1,\ldots,n\} 
\text{ for some cell } Q \text{ of } A_\omega\right. \right\}.\]
In particular, by our smoothness assumption on $\pi^{-1}(0)\cap\tZ$, 
Corollary \ref{cor:simplex} tells us that $\pi^{-1}(0)\cap\tZ$ is actually\\
\scalebox{.85}[1]{$\displaystyle{\left\{[p_1:\cdots:p_{2N}]\!\in\!Y_{\tilde{A}} \; \; 
\left| \; \;  \sum\limits_{a_j\in Q} c_{i,{a_j}}p_j\!=\!0 \text{ for all } 
i\!\in\!\{1,\ldots,n\} \text{ for some full-dimensional cell } Q \text{ of } A_\omega\right. 
\right\}}$.}
Since $A_\omega$ is a triangulation, Corollary \ref{cor:simplex} and Lemma 
\ref{lemma:veronese} (along with another application Hermite factorization) 
tells us that 
\[\#\pi^{-1}(0)\cap\tZ\!=\!\sum\limits_{Q \text{ a full-dimensional 
cell of } A_\omega} \frac{\vol(Q)}{h}\!=\!\frac{\vol(A)}{h},\] so we are 
done. \qed 

\begin{rem} 
Our proof generalizes quite easily to arbitrary algebraically closed fields 
and positive characterstic, e.g., the algebraic closure of a finite field. 
One need only use a little algebra to extend Lemma \ref{lemma:curve} to 
algebraically closed fields and then one can use the same proof above 
almost verbatim. \dia 
\end{rem} 

\subsection{Path Following, Compactness, and the Degenerate Case of Kushnirenko}\mbox{}\\
Let us now allow degeneracies for the zero set of $F$ and arrive at a strengthening of 
Kushnirenko's Theorem. 
\begin{thm} 
\label{thm:hardkush}  
Following the notation of Theorem \ref{thm:kush}, let $Z_A$ be the zero set of 
$F$ in $Y_A$, and let $\{Z_i\}$ be the collection of path-connected components 
of $Z_A$. Then there is a natural, well-defined positive {\bf intersection multiplicity} 
$\mu : \{Z_i\} \longrightarrow \N$ such that $\sum_i\mu(Z_i)\!=\!\vol(A)$ 
and $\mu(Z_i)\!=\!1$ if $Z_i$ is a non-degenerate root.  
\end{thm} 

We actually have all the technical preliminaries we'll need, except for one 
last simple proposition on path-connectedness. 
\begin{prop}
\label{prop:path}
If $\cH$ is any algebraic hypersurface in $\C^N$ then 
$B\!\setminus\!\cH$ is path-connected for any open ball $B$ 
in $\C^N$. \qed 
\end{prop}

\noindent
{\bf Proof of Theorem \ref{thm:hardkush}:} Let $N\!:=\!\#A$ as usual and 
let $P(A)$ denote the space of all polynomials in $\C[x_1,\ldots,x_n]$ with support contained 
in $A$. Clearly, we can naturally identify $P(A)$ with $\C^N$. Since zero sets of polynomials 
are unchanged under scaling of the coefficients, we will then let 
$\Pro(\underset{n}{\underbrace{A,\ldots,A}})\!:=\!(\Pro^{N-1}_\C)^n$ be the 
space we'll use to consider our possible $F$. 

Note now that if $F\!\in\!(\Pro^{N-1}_\C)^n\!\setminus\!
\Delta_{(\underset{n}{\underbrace{A,\ldots,A}})}$ then we are done by 
Theorem \ref{thm:smoothkush} (simply setting $\mu(Z_i)\!=\!1$ for every 
root $Z_i$). Indeed, since 
$(\Pro^{N-1}_\C)^n\!\setminus\!
\Delta_{(\underset{n}{\underbrace{A,\ldots,A}})}$ is path-connected by 
Proposition \ref{prop:path}, the Implicit Function Theorem tells us that 
$F$ had better have the same number of roots in $Y_A$ as an 
$F$ which has smooth zero set and no roots at toric infinity. 

Essentially the same idea can be used for 
$F\!\in\!\Delta_{(\underset{n}{\underbrace{A,\ldots,A}})}$. In particular, 
for such $F$, let $F^{(i)}$ be any sequence such that 
$F^{(i)} \longrightarrow F$. Then, letting $Z^{(i)}$ be the zero set of 
$F^{(i)}$ in $Y_A$, define $\zeta$ to be any limit point of $\{Z^{(i)}\}$. 
We must then have $F(\zeta)\!=\!0$ and thus $Z_A$ must be non-empty. 

Now let $\{U_i\}$ be disjoint open sets with $Z_i\!\subset\!U_i$ for all 
$i$. (Such a set of open sets must exist since $Y_A$ is compact and 
the $Z_i$ must be of positive distance from each other, using the usual 
distance in $\Pro^{N-1}_\C$.) Note then that $Y_A\!\setminus\!\bigcup_i U_i$ 
must be compact. By the continuity of $F$ as a function 
of its variables {\bf and} coefficients, there 
must then be a ball $B$ about $F$ in $(\Pro^{N-1}_\C)^n$ 
such that the roots of any $G$ are contained in $\bigcup_i U_i$. 

We may now define $\mu(Z_i)$ as follows: Take any $G\!\in\!B\setminus 
\Delta_{(\underset{n}{\underbrace{A,\ldots,A}})}$ and define 
$\mu(Z_i)$ to be the number of roots of $G$ in $U_i$. Since $B\setminus
\Delta_{(\underset{n}{\underbrace{A,\ldots,A}})}$ is path-connected by 
Proposition \ref{prop:path}, the Implicit Function Theorem tells us that 
the number of roots is independent of whatever $G\!\in\!B\setminus
\Delta_{(\underset{n}{\underbrace{A,\ldots,A}})}$ we picked. 

\begin{rem} 
Essentially the same argument was used Mike Shub in \cite{shub} to prove an extended version 
of B\'ezout's Theorem (see also \cite[Pg.\ 199]{bcss}). Note that neither theorem 
generalizes the other, but the theorems overlap in the special case 
where $A$ is the set of integral points in a scaled standard simplex. 
On the other hand, our next theorem will generalize both Kushnirenko's Theorem 
and B\'ezout's Theorem.  The elegance of Shub's approach is that it gives a rigourous and 
simple approach to intersection theory for a broad class of polynomial systems. 
\dia 
\end{rem} 

\begin{rem} 
Combining Theorems \ref{thm:hardkush} and \ref{thm:smoothkush} we immediately 
obtain our earlier simpler statement of Kushnirenko's Theorem (Theorem \ref{thm:kush}). \dia 
\end{rem} 

\section{Multilinearity and Reducing Bernstein's Theorem to 
\mbox{Kushnirenko's} Theorem} 
\label{sec:multi} 
The big question now is how to count the roots of a {\bf mixed} 
polynomial system, since being unmixed is such a strong restriction. 
Toward this end, let us consider another consequence of the 
basic properties of discriminant varieties. 
\begin{lemma} 
\label{lemma:add} 
Let $F$ and $G$ be any $n\times n$ polynomial systems with 
support contained in $(A_1,\ldots,A_n)$ component-wise. Then, generically, $F$ and $G$ share 
no roots in $\Csn$. Furthermore, the number of roots of 
$F$ is generically a fixed constant. \qed 
\end{lemma} 

As an immediate consequence, we obtain the following preliminary 
answer to our big question. 
\begin{dfn} 
Let $S_1,\ldots,S_k$ be any subsets of $\Rn$. Then 
their {\bf Minkowski sum} is simply $S_1+\cdots+S_k\!:=\!
\{y_1+\cdots+y_k \; | \; y_i\!\in\!S_i \text{ for all } i\}$.  \dia 
\end{dfn} 
\noindent
It is easily proved that $\newt(fg)\!=\!\newt(f)+\newt(g)$ (once one 
observes that the vertices of $\newt(fg)$ are themselves Minkowski 
sums of vertices of $\newt(f)$ and $\newt(g)$). So it should 
come as no surprise that Minkowksi sums will figure importantly in our 
discussion relating polyhedra and polynomials. 
\begin{lemma}
\label{lemma:multi} 
Let $\cN(A_1,\ldots,A_n)$ denote the generic number of roots in $\Csn$ of an 
$n\times n$ polynomial system $F$ with 
support $(A_1,\ldots,A_n)$. Then $\cN(A_1,\ldots,A_n)$ is a  
non-negative symmetric function of $\conv(A_1),\ldots,\conv(A_n)$ 
which is multilinear with respect to Minkowski sum. 
\end{lemma} 

\noindent
{\bf Proof:} That $\cN(A_1,\ldots,A_n)$ is a well-defined non-negative symmetric function of 
$A_1,\ldots,A_n$ is clear (thanks in part to the last part of 
Lemma \ref{lemma:add}). The formula for $\cN(A_1,\ldots,A_n)$ in the unmixed case then follows 
immediately from Theorem \ref{thm:smoothkush}. Translation invariance 
follows easily since the roots of $F$ in $\Csn$ are the same as the roots of 
$(x^{a_1}f_1,\ldots,x^{a_n}f_n)$ in $\Csn$. Defining $x^[u_{ij}]\!:=\!(
x^{u_{11}}_1\cdots x^{u_{n1}},\ldots,x^{u_{1n}}_1\cdots x^{u_{nn}})$ for any 
$n\times n$ matrix $[u_{ij}]$, it is then easily checked that 
$\supp(F(x^U))\!=\!(UA_1,\ldots,UA_n)$ and $U$ 
{\bf unimodular} (cf.\ Definition \ref{dfn:hermite}) $\Longrightarrow$ the map 
$x\mapsto x^U$ is an analytic bijection of $\Csn$ into itself. 
So it is clear that $\cN(UA_1,\ldots,UA_n)\!=\!\cN(A_1,\ldots,A_n)$. 

We thus need only show that $\cN$ is a multilinear function of the convex hulls. 
To see the multilinearity, note that the zero set of $(f_1\bar{f}_1,f_2,\ldots,f_n)$ in 
$\Csn$ is exactly the union of the zero sets of $(f_1,f_2,\ldots,f_n)$ 
$(\bar{f}_1,f_2,\ldots,f_n)$. So by the first part of Lemma \ref{lemma:add}, and the symmetry 
of $\cN$, multlinearity follows. 

Recall now the {\bf polarization identity}: 
\[ m(x_1,\ldots,x_n)\!=\!\sum\limits_{\emptyset\neq I\subseteq\{1,\ldots,n\}}
(-1)^{n-\#I} \begin{pmatrix} n \\ \#I\end{pmatrix} 
m\left(\sum\limits_{i\in I} x_i,\ldots,\sum\limits_{i\in I} x_i\right), \] 
valid for any symmetric multilinear function. (The identity is 
not hard to prove via inclusion-exclusion \cite{gkp}. See also \cite{goldman} 
in this volume for another point of view.) Therefore, we must have 
\[ \cN(A_1,\ldots,A_n)\!=\!\sum\limits_{\emptyset\neq I\subseteq\{1,\ldots,n\}}
(-1)^{n-\#I} \begin{pmatrix} n \\ \#I \end{pmatrix} 
\cN\left(\sum\limits_{i\in I} A_i,\ldots,\sum\limits_{i\in I} A_i\right), \] 
and thus $\cN(A_1,\ldots,A_n)$ depends only the convex hulls of $A_1,\ldots,A_n$. \qed 

So we have answered our big question, assuming we know a function $\cM(P_1,\ldots,P_n)$,  
defined on $n$-tuples $(P_1,\ldots,P_n)$ of polytopes in $\Rn$, that satisfies the obvious 
analogues of the properties of $\cN(A_1,\ldots,A_n)$ specified in Lemma \ref{lemma:multi}. 
However, such a function indeed exists: it is called the {\bf mixed volume} and 
we denote it by $\cM(\cdot)$. Abusing notation slightly by setting 
$\cM(A_1,\ldots,A_n)\!:=\!\cM(\conv(A_1),\ldots,\conv(A_n))$, we immediately obtain 
the following result. 
\begin{thm}
[Bernstein's Theorem] 
\label{thm:bernie} 
Suppose $F$ is any $n\times n$ polynomial system with fixed support $A_1,\ldots,A_n$.  
Then $F$ generically has exactly $\cM(A_1,\ldots,A_n)$ roots in $\Csn$. 
\end{thm} 

Of course, we now appear to have an even bigger question: what is mixed volume? 
This we now answer. 

{\sc Remark 4.0.1}

\vspace{-.5cm}
\noindent
\begin{minipage}[t]{2in}
\begin{picture}(40,30)(0,0)
\put(0,-90){\epsfig{file=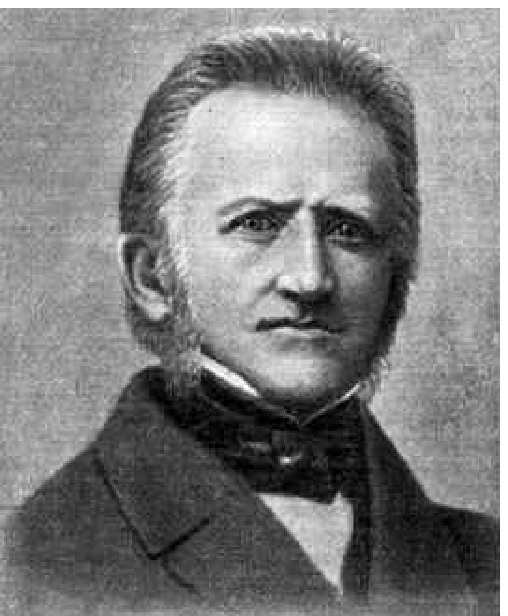,height=1.4in}\hspace{.2cm}
\epsfig{file=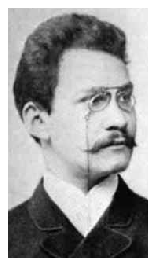,height=1.4in}}
\put(15,-100){{\small Felix Minding}}
\put(100,-100){{\small \scalebox{1}[.7]{$\begin{matrix}\text{Hermann} \\ \text{Minkowski}
\end{matrix}$}}}
\end{picture} 
\end{minipage}
\hspace{.1cm}
\begin{minipage}[t]{2.5in}
{\em 
David N.\ Bernstein was the first to prove Theorem \ref{thm:bernie} in a slightly 
stronger form, along with an algebraic condition for when $F$ would have 
exactly mixed volume many roots \cite{bernie}. Interestingly, Felix Minding appears to have 
been the first to prove the special case $n\!=\!2$ in 1841, and mixed volume wasn't 
even defined until near the end of the $19\thth$ century by Hermann Minkowski. \dia} 
\end{minipage} 

\medskip
\noindent
We also point out that Minkowski was born on 22 June, 1864, in a town named Alexotas. 
This town used to belong to what was the Russian empire at the time but is now 
the Lithuanian city of Kaunas. 
 
\section{Mixed Subdivisions and Mixed Volumes {}from Scratch} 
\label{sec:mixed} 

Let us begin with an illustration of one of the simplest non-trivial examples of a 
Minkowski sum:\\ 
\raisebox{-.4cm}{\epsfig{file=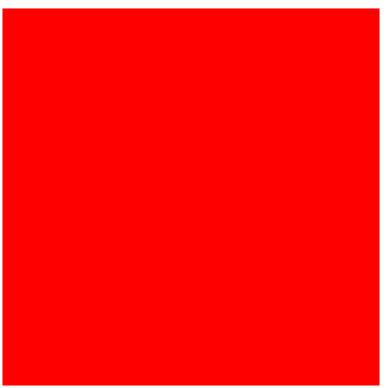,height=.5in}}\hspace{.1cm}
\scalebox{1.4}[1.4]
{{\huge $\pmb{+}$}}
\raisebox{-.4cm}{\epsfig{file=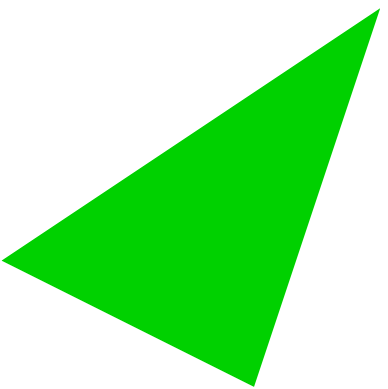,height=.5in}}
\scalebox{1.4}[1.4]
{{\huge $ \ \pmb{=} \ $}}
\raisebox{-1cm}{\epsfig{file=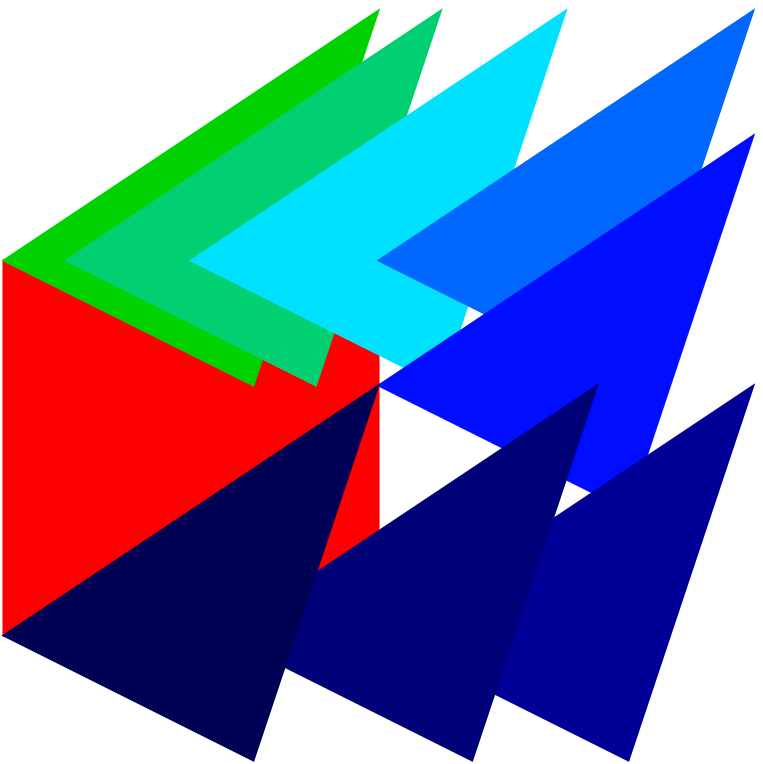,height=1in}}
\scalebox{1.4}[1.4]
{{\huge $ \ \pmb{=} \ $}}
\raisebox{-1cm}{\epsfig{file=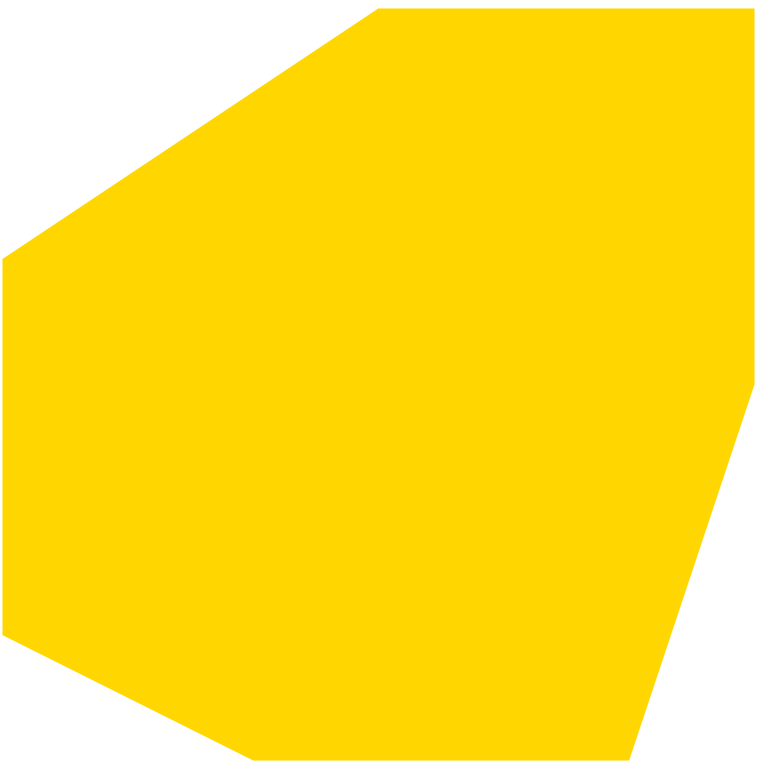,height=1in}}

There are many different defintions of mixed volume but the two most important 
use Minkowski sums in an essential way. More to the point, if one can subdivide 
$P_1+\cdots+P_n$ in a special way, then one is well on the way to computing 
mixed volume. This is where {\bf mixed subdivisions} enter.  
\begin{dfn} \cite{hs} 
Given polytopes $P_1,\ldots,P_k\!\subset\!\Rn$, 
a {\bf subdivision of $\pmb{(P_1,\ldots,P_k)}$} is a finite collection of $k$-tuples
$\{(C^\alpha_1,\ldots,C^\alpha_k)\}_{\alpha\in S}$ 
satisfying the following axioms:
\begin{enumerate}
\item{$\bigcup_{\alpha\in S} 
C^\alpha_i\!=\!P_i$ for all $i$} 
\item{$C^\alpha_i\cap C^\beta_i$
is a face of both $C^\alpha_i$ and 
$C^\beta_i$ for all $\alpha,\beta,i$. }
\item{\scalebox{.9}[1]
{$C^\beta_i$ a face of $C^\alpha_i$ for all $i \Longrightarrow
C^\beta_1,\ldots,C^\beta_n$ have a common inner normal.}}
\end{enumerate}
Furthemore, if we have in addition that 
$\sum_i \dim C^\alpha_i\!=\!\dim \sum_i C^\alpha_i$
for all $\alpha$, then we call $\{(C^\alpha_1,\ldots,C^\alpha_k)\}_{\alpha\in S}$  
a {\bf mixed subdivision}. \dia 
\end{dfn}
\begin{ex}\mbox{}\\ 
\scalebox{.9}[.9]{\begin{picture}(300,200)(0,20)
\put(30,0){\epsfig{file=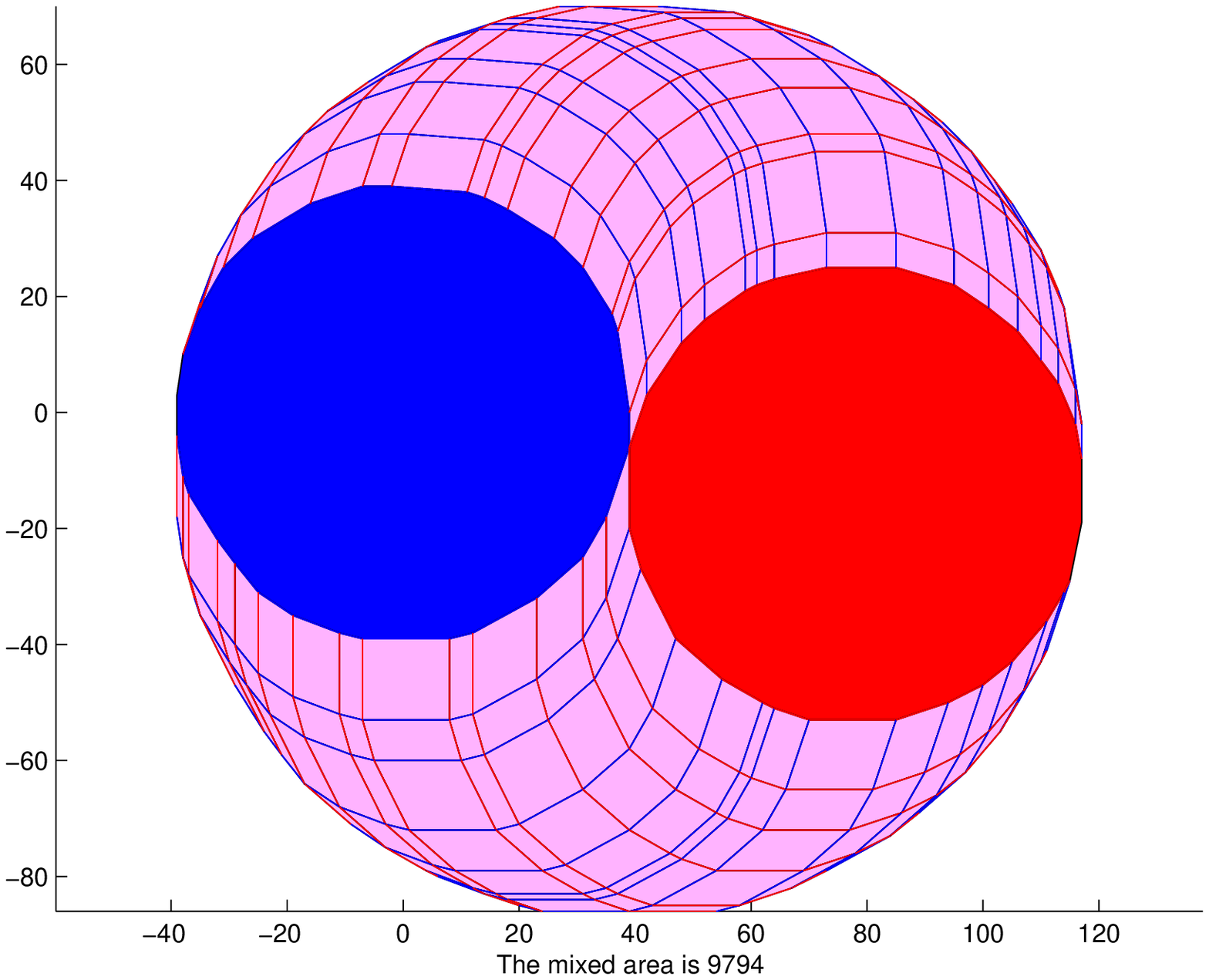,height=3in}}
\put(90,125){\textcolor{yell}{$P_1+v_2$}}
\put(80,110){{\normalsize \textcolor{yell}{$(P_1,v_2)$ {\bf unmixed}}}}
\put(158,118){\textcolor{yell}{{\normalsize $v_1$}}}
\put(166,116){\textcolor{yell}{{\tiny $\bullet$}}}
\put(169,114){\textcolor{yell}{{\normalsize $v_2$}}}
\put(195,110){\textcolor{yell}{$v_1+P_2$}}
\put(180,95){{\normalsize \textcolor{yell}{$(v_1,P_2)$ {\bf unmixed}}}}
\put(220,170){\begin{rotate}{45}
\scalebox{1.5}[1]{$\pmb{\longleftarrow}$}
\end{rotate}}
\put(245,200){{\large 
$\mathrm{Edge}_1+\mathrm{Edge}_2$}}
\put(248,184){\scalebox{.8}[1]{{\normalsize 
$(\mathrm{Edge}_1,\mathrm{Edge}_2)$ is {\bf mixed}}}}
\put(258,170){\scalebox{.8}[1]{{\normalsize 
...or of {\bf type} $(1,1)$}}}
\end{picture}}

\medskip
\bigskip
Here we see a very special kind of subdivision $\{Q_i\}$ of the Minkowski sum 
of two polygons $P_1$ and $P_2$, each  
with many vertices. In particular, the subdivision of $P_1+P_2$ above is built 
in such a way as to encode a {\bf mixed} subdivision $\{(C^{\alpha}_1,C^{\alpha}_2)$ 
of $(P_1,P_2)$. In particular, we see that each $P_i$ has a distinguished vertex $v_i$, and 
that we can read off a mixed subdivision of $(P_1,P_2)$ as follows: there are 
two cells $(P_1,v_2)$ and $(v_1,P_2)$, corresponding to the two cells  
$P_1+v_2$ and $v_1+P_2$ of $\{Q_i\}$. The remaining cells of $\{(C^{\alpha}_1,C^{\alpha}_2)$ 
are parallelograms of the form $(E_1,E_2)$ where $E_i$ is an edge of 
$P_i$ for all $i$. \dia 
\end{ex} 

It is easily verified that any subdivision of $(P_1,\ldots,P_k)$ immediately 
induces a subdivision of $(\lambda P_1,\ldots,\lambda P_k)$, for any 
$\lambda_1,\ldots,\lambda_k\!\geq\!0$. 
\begin{ex}\mbox{}\\
\scalebox{.9}[.9]{
\begin{picture}(300,200)(0,20)
\put(30,0){\epsfig{file=shockwave.eps,height=3in}}
\put(70,122){\scalebox{.85}[1]{\textcolor{yell}{{\normalsize
$\area(\lambda P_1)\!=\!\lambda^2\area(P_1)$}}}}
\put(170,104){\scalebox{.85}[1]{\textcolor{yell}{{\normalsize
$\area(\mu P_2)\!=\!\mu^2\area(P_2)$}}}}
\put(220,170){\begin{rotate}{45}
\scalebox{1.5}[1]{$\pmb{\longleftarrow}$}
\end{rotate}}
\put(245,200){\scalebox{.65}[1]{{\normalsize 
$\area(\text{Mixed Cell of } (\lambda P_1,\mu P_2))$}}}
\put(248,184){\scalebox{.65}[1]{{\normalsize 
$=\!\lambda\mu\area(\text{Mixed Cell of } (P_1,P_2))$
}}}
\end{picture}}

\medskip
\bigskip
Note in particular that the areas of the cells of our induced subdivision  
of $\lambda P_1+\mu P_2$ scale according to their {\bf type}. In particular, 
it is clear that for our above example, $\area(\lambda P_1+\mu P_2)
\!=\!\area(P_1)\lambda^2+M\lambda\mu+\area(P_2)\mu^2$, where $M$ is the 
sum of the areas of all the {\bf mixed} cells (the parallelograms). \dia 
\end{ex} 

\begin{dfn} 
Following the notation above, the {\bf type} of a cell $(C^{\alpha}_1,\ldots,C^{\alpha}_k)$ 
of a subdivision of $(P_1,\ldots,P_n)$ is simply the vector 
$(\dim C^{\alpha}_1,\ldots,\dim C^{\alpha}_k)$. In particular, the 
cells of type $(1,\ldots,1)$ are called {\bf mixed cells}. 
\dia 
\end{dfn} 

The following two lemmata are then immediate. The first follows from a slight 
modification of the proof of Lemma \ref{lemma:lift}, while the second follows almost 
immediately from the first. 
\begin{lemma}
\label{lemma:msd}
Following the notation of Definition \ref{dfn:lift}, recall that 
\[ \pi : \R^{n+1} \longrightarrow \Rn\] is the natural projection which 
forgets the last coordinate. 
Then, given finite point sets $A_1,\ldots,A_n\!\subset\!\Zn$ and lifting functions 
$\omega_i$ for $A_i$ for all $i$, the collection 
$(A_1,\ldots,A_n)_\omega\!:=\!\{(\pi(\conv(\hat{A}_1)^w),\ldots,\pi(\conv(\hat{A}_n)^w)) 
\; | \; w\!\in\!\Rn\!\setminus\!\{\bO\}\}$ always forms a subdivision 
of $(\conv(A_1),\ldots,\conv(A_n))$ --- the {\bf subdivision of 
$\pmb{(\conv(A_1),\ldots,\conv(A_n))}$ induced by $\pmb{(\omega_1,\ldots,\omega_n)}$}.
In particular, for fixed $(A_1,\ldots,A_n)$, $(A_1,\ldots,A_n)_\omega$ will 
generically be a {\bf mixed} subdivision. \qed 
\end{lemma} 
\begin{lemma} 
\label{lemma:vol} 
For $\lambda_1\ldots,\lambda_n\!\geq\!0$, and any polytopes 
$P_1,\ldots,P_n\!\subset\!\Rn$, the quantity  
$Q(\lambda_1,\ldots,\lambda_n):=\vol\left(\sum^n_{i=1}\lambda_i P_i\right)$
is a homogeneous polynomial of degree $n$ with nonnegative coefficients. \qed 
\end{lemma} 

We then at last arrive at the following definition of the mixed volume. 
\begin{dfn} 
Given any polytopes $P_1,\ldots,P_n\!\subset\!\Rn$, their mixed volume 
is the coefficient of $\lambda_1\lambda_2\cdots \lambda_n$ in the 
above polynomial $Q(\lambda_1,\ldots,\lambda_n)$. \dia 
\end{dfn} 
\begin{ex}[The Unmixed Case]
It is easily checked that $\cM(P,\ldots,P)\!=\!\vol(P)$. Note also 
that the multilinearity of $\cM(\cdot)$ with respect to Minkowski 
sum also follows immediately from the preceding definition. \dia 
\end{ex} 
\begin{ex}[Line Segments]
It is also easily checked that $\cM(\{0,a_1\},\ldots,\{0,a_n\})\!=\!|\det[
a_1,\ldots,a_n]|$, where $a_1,\ldots,a_n$ are any points in $\Rn$ and 
$[a_1,\ldots,a_n]$ is the matrix whose columns are $a_1,\ldots,a_n$. \dia 
\end{ex} 

The next two characterizations follow easily from the last two lemmata, and 
inclusion-exclusion \cite{gkp}. 
\begin{lemma} 
For any mixed subdivision $\{(C^{\alpha}_1,\ldots,C^{\alpha}_n)\}$ 
of $(P_1,\ldots,P_n)$, 
\[\cM(P_1,\ldots,P_n):=\!\!\!\!\!\!\!\!\sum
\limits_{\substack{ (C_1,\ldots,C_n) \\
\text{ a cell of type } (1,\ldots,1) }} \!\!\!\!\!\!\!\!
\vol\left(\sum_i C_i\right).\] 
Furthermore, we have $\displaystyle{\cM(P_1,\ldots,P_n):=\!\!\!\!\!\!\sum_{\emptyset\neq 
I\subseteq \{1,\ldots,n\}}  (-1)^{n-\#I}\vol\left(\sum_{i\in I} P_i\right)}$. \text{\dia} 
\end{lemma} 

\begin{ex}[Cornered Spikes] 
A less trivial puzzle is the following formula: 
$\cM(\{\bO,a_{11}e_1,a_{1n}e_n\},\ldots,\{\bO,a_{11}e_1,a_{1n}e_n\})\!=\!\max\limits_\sigma
\left\{\prod\limits^n_{i=1} a_{i\sigma(i)}\right\}$, where $a_{ij}$ are any 
non-negative real numbers and $\sigma$ ranges over 
all permutations of $\{1,\ldots,n\}$. This gives some indication that 
mixed volume includes many simple functions as a special case. 
However, our next example shows that mixed volume includes 
rather non-trivial functions as well. \dia 
\end{ex} 
\begin{ex}[Bricks, a.k.a.\ the fine multigraded case] 
Via multilinearity, it easily follows that 
$\cM([0,d_{11}]\times \cdots \times [0,d_{1n}],
\ldots,[0,d_{n1}]\times \cdots \times [0,d_{nn}])\!=\!\perm[d_{ij}]$, 
where $\perm$ denotes the {\bf permanent}.\footnote{Recall that this 
function can be defined as the variant of the determinant where all 
alternating signs in the full determinant expansion are replaced by 
$+1$'s.} In particular, this immediately shows that computing 
mixed volume is $\#\mathbf{P}$-hard \cite{papa,dgh}. \dia 
\end{ex} 

Let us now finally prove Theorem 2.\\
{\bf Proof of Theorem 2:} Note that by Bernstein's Theorem, 
it suffices to find an algorithm for computing $\cM(A_1,A_2)$ with 
bit complexity $\cO(b\bar{N}+\bar{N}\log \bar{N})$. 
The main idea of the proof can then already be 
visualized in the first mixed subdivision we illustrated: 
one computes the mixed area of $(A_1,A_2)$ by first efficiently computing the 
convex hulls of $A_1$ and $A_2$, and then expressing the sum of the 
areas of the mixed cells compactly {\bf without} building the entire 
mixed subdivision. This is not a contradiction, provided one views the mixed cells in the 
right way.

More precisely, first recall that the convex hulls of $A_1$ and $A_2$ can 
be computed within $O(\bar{N}\log \bar{N})$ bit operations, via the usual 
well-known $2$-dimensional convex hull algorithms \cite{preparata}. In particular, 
with this much work, we can already assume we know the inner edge normals 
of $P_1\!:=\!\conv(A_1)$ and $P_2\!:=\!\conv(A_2)$, and the vertices of 
$P_1$ and $P_2$ in counter-clockwise order. 

Let us then pick a vertices $v_1\!\in\!P_1$ and $v_2\!\in\!P_2$ such that 
their angle cones are disjoint. Then there is a mixed subdivison 
(which we will never calculate explicitly!) with exactly $2$ 
non-mixed cells --- $(P_1,v_2)$ and $(v_1,P_2)$ --- and several other 
mixed cells. (This is easily seen by picking a lifting function $\omega_1$ for 
$P_1$ that is identically zero, and a linear lifting function $\omega_2$ for $P_2$ 
that is minimized at $v_2$ and is constant on a line that intersects 
the angle cones of $v_1$ and $v_2$ only at the origin.) 

Note then that the union of the mixed cells can be partition into a 
union of strips. In particular, by construction, there are 
disjoint contiguous sequences of edges $(E^{(i)}_1,\ldots,E^{(i)}_{a_i})$ and 
$(E^{(i')}_1,\ldots,E^{(i')}_{a_{i'}})$, with $E^{(i)}_1$ and $E^{(i')}_1$ 
incident to $v_i$, for all $i$. Furthermore, every mixed cell of 
$(A_1,A_2)_\omega$ is of the form $(E^{(1)}_i,E^{(2)}_j)$ or 
$(E^{(1')}_i,E^{(2')}_j)$, and every $E^{(i)}_j$ and $E^{(i')}_j$ is 
incident to some mixed cell. 

The partition into strips then arises as follows: the mixed cells of 
$(A_1,A_2)_\omega$ can be partitioned into lists of one of the following two forms: 
\[ (E^{(1)}_j,E^{(2)}_{m_j}),\ldots,(E^{(1)}_1,E^{(2)}_{n_j} ) \]
\[ (E^{(1')}_j,E^{(2')}_{m_j}),\ldots,(E^{(1')}_1,E^{(2')}_{n_j} ) \] 
where $j\!\in\!\{1,\ldots,a_1\}$ (resp.\ $j\!\in\!\{1,\ldots,a_{1'}\}$),  
$m_j\!\leq\!n_j$, and $n_j\!\leq\!a_2$ (resp.\ $n_j\!\leq\!a_{2'}$. 
In particular, the union of the mixed cells in any such list is 
simply the Minkowski sum of a continuous portion of the boundary 
of $P_2$ and an edge of $P_1$, and its area can thus be expressed as the 
absolute value of a determinant of differences of vertices of the $P_i$. 
Furthermore, each formula can easily be found by a binary search on the 
sorted edge normals using $O(\bar{N}\log\bar{N})$ comparisons. 

Since there are no more than $\bar{N}$ such strips, the total work we do 
is bounded above by the specified complexity bound, so our upper bound 
is proved. 

To obtain our lower bound, note that the mixed area of $(A_1,A_2)$ is 
zero iff [[$P_1$ or $P_2$ is a point] or [$P_1$ and $P_2$ are parallel 
line segments]. So just knowing whether the mixed area is 
positive or not amounts to a rank computation on a matrix 
of size $O(\bar{N})$ and thus can take no less than $\Omega(bN)$ 
bit operations in the worst case \cite{bcs}. \qed  

\section{A Stronger Bernstein Theorem Via Mixed Subdivisions} 
\label{sec:direct} 

Here we prove the following generalization of Theorem \ref{thm:hardkush}. 
It is at this point that we will use a slightly more high-brow type of 
toric variety: the toric variety $X_P$ corresponding to a polytope $P$. 
In essence, the key properties that we needed from $Y_A$ (that it 
compactify $\Csn$ and have a partition into orbits corresonding to 
the faces of a polytope) continue to hold for $X_P$. We make 
this change to avoid technicalities in defining the zero set of 
$F$ in $Y_{A_1+\cdots+A_n}$ when $F$ is mixed. Since $X_P$ is 
discussed elsewhere in this volume at greater length \cite{cox,frank}, we 
proceed with the statement of our theorem. 
\begin{thm}
\label{thm:hardbernie}
Following the notation of Theorem \ref{thm:bernie}, let $P\!:=\!\conv(A_1)+\cdots+\conv(A_n)$,  
let $Z_P$ be the zero set of $F$ in $X_P$, and let $\{Z_i\}$ be the collection of 
path-connected components of $Z_P$. Then there is a natural, well-defined positive {\bf 
intersection multiplicity} $\mu : \{Z_i\} \longrightarrow \N$ such that 
$\sum_i\mu(Z_i)\!=\!\vol(A)$ and $\mu(Z_i)\!=\!1$ if $Z_i$ is a non-degenerate root. 
\end{thm}

The proof will be almost exactly the same as that of our extended version 
of Kushnirenko's Theorem, so let us first see an illustration of a toric 
deformation for a {\bf mixed} system. 

\begin{ex} 
Take $n\!=\!2$ and
\[ f_1(x,y)\!:=\!c_{1,\bO}+c_{1,(\alpha,0)}x^\alpha+c_{1,(0,\beta)}y^\beta
+c_{1,(\alpha,\beta)}x^\alpha y^\beta \]
\[ f_2(x,y)\!:=\!c_{2,\bO}+c_{2,(\gamma,0)}x^\gamma+c_{2,(0,\delta)}y^\delta
+c_{2,(\gamma,\delta)}x^\gamma y^\delta. \] 
By Bernstein's Theorem, the number of roots should be $\alpha\delta+\beta\gamma$, 
so let us try to prove this.

Let us take the following lifting of $F$: 
\[ \hat{f}_1(x,y,t)\!:=\!c_{1,\bO}+c_{1,(\alpha,0)}x^\alpha t
+c_{1,(0,\beta)}y^\beta t 
+c_{1,(\alpha,\beta)}x^\alpha y^\beta \]
\[ \hat{f}_{2}(x,y)\!:=\!c_{2,\bO}t+c_{2,(\gamma,0)}x^\gamma
+c_{2,(0,\delta)} y^\delta +c_{2,(\gamma,\delta)}x^\gamma
y^\delta t \] 
In particular, we see that there will be exactly one mixed cell for 
$(A_1,A_2)_{\omega}$ and its corresponding initial term system will be 
\[ \init_{(0,0,1)}(\hat{F})(x,y,t)=(c_{1,\bO}+
c_{1,(\alpha,\beta)}x^\alpha y^\beta,c_{2,(\gamma,0)}x^\gamma+c_{2,(0,\delta)}
y^\delta) \] 
The lifted Newton polytopes and induced subdivisions appear below.\\ 
\begin{picture}(300,200)(0,20)
\put(0,0){\epsfig{file=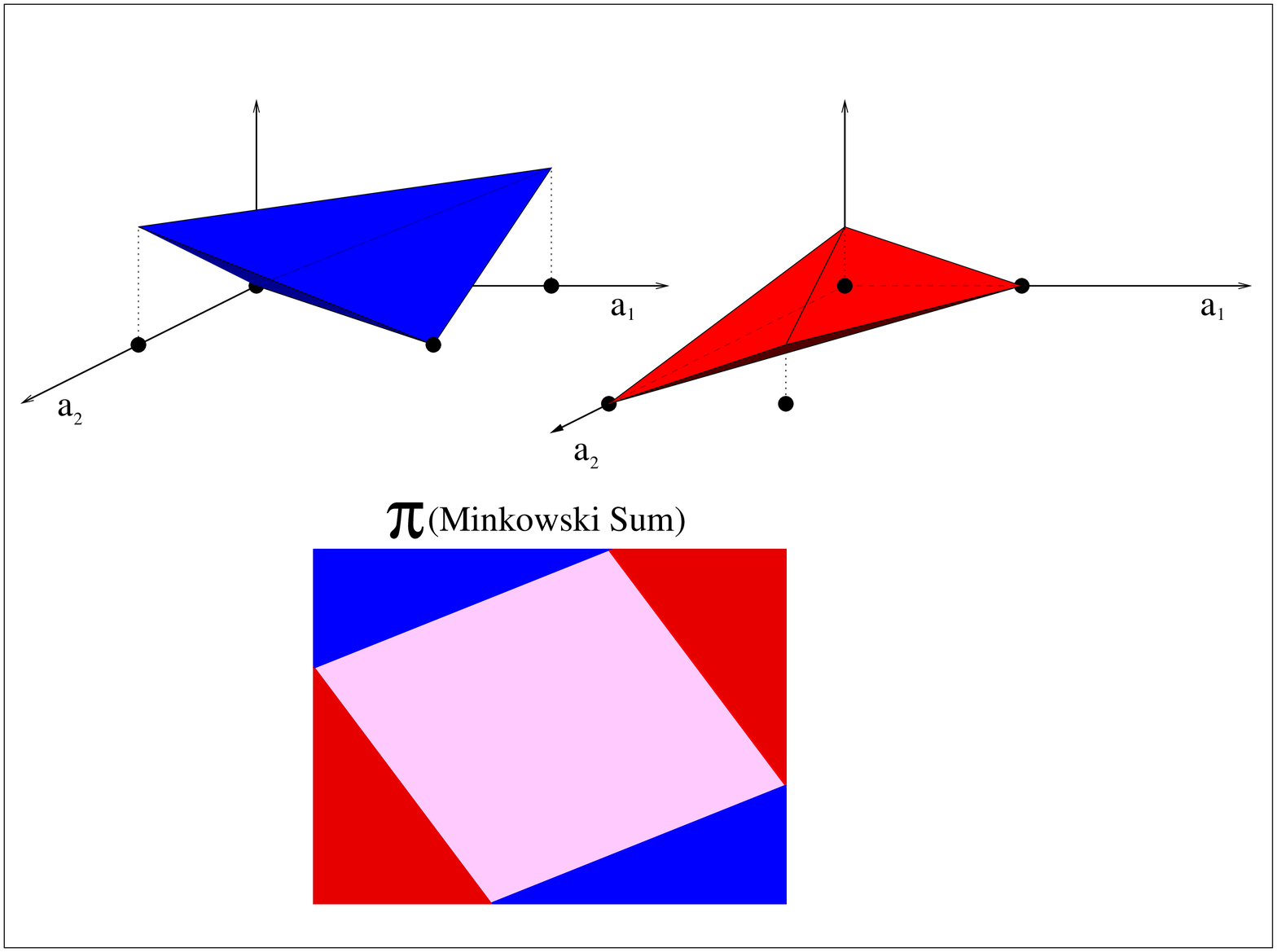,height=3.3in}}
\end{picture} 

\bigskip
The idea of our proof of Bernstein's Theorem then mimics our earlier proof of Kushnirenko's 
Theorem: our lifting induces a lifted version $\hat{P}\!=\!\conv(\hat{A}_1)+
\conv(\hat{A}_2)$ of $P$ and we'll then try to build a map from our 
lifted zero set to the projective line. To do so, we'll define $\tilde{P}\!:=\!\hat{P}\times 
[0,1]$ and this is illustrated below.\\ 
\begin{picture}(300,190)(0,16)
\put(0,40){\scalebox{.5}[1]{\epsfig{file=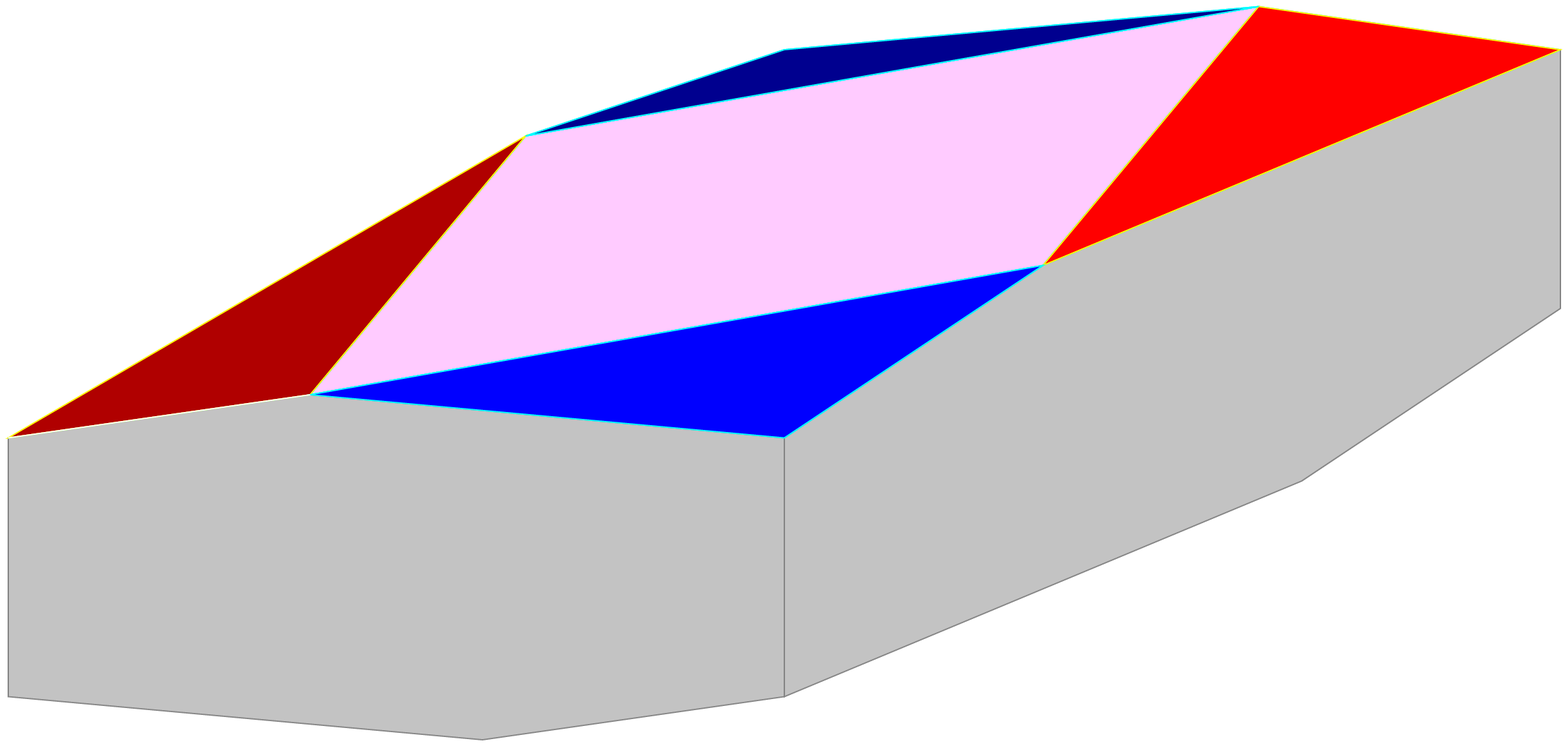,height=2in}}}
\put(18,72){\scalebox{.8}[1]{{\LARGE $X_{\hat{P}+I}$}}}
\put(140,80){{\LARGE $\pmb{\approx}$}}
\put(170,-20){\scalebox{.5}[1]{\epsfig{file=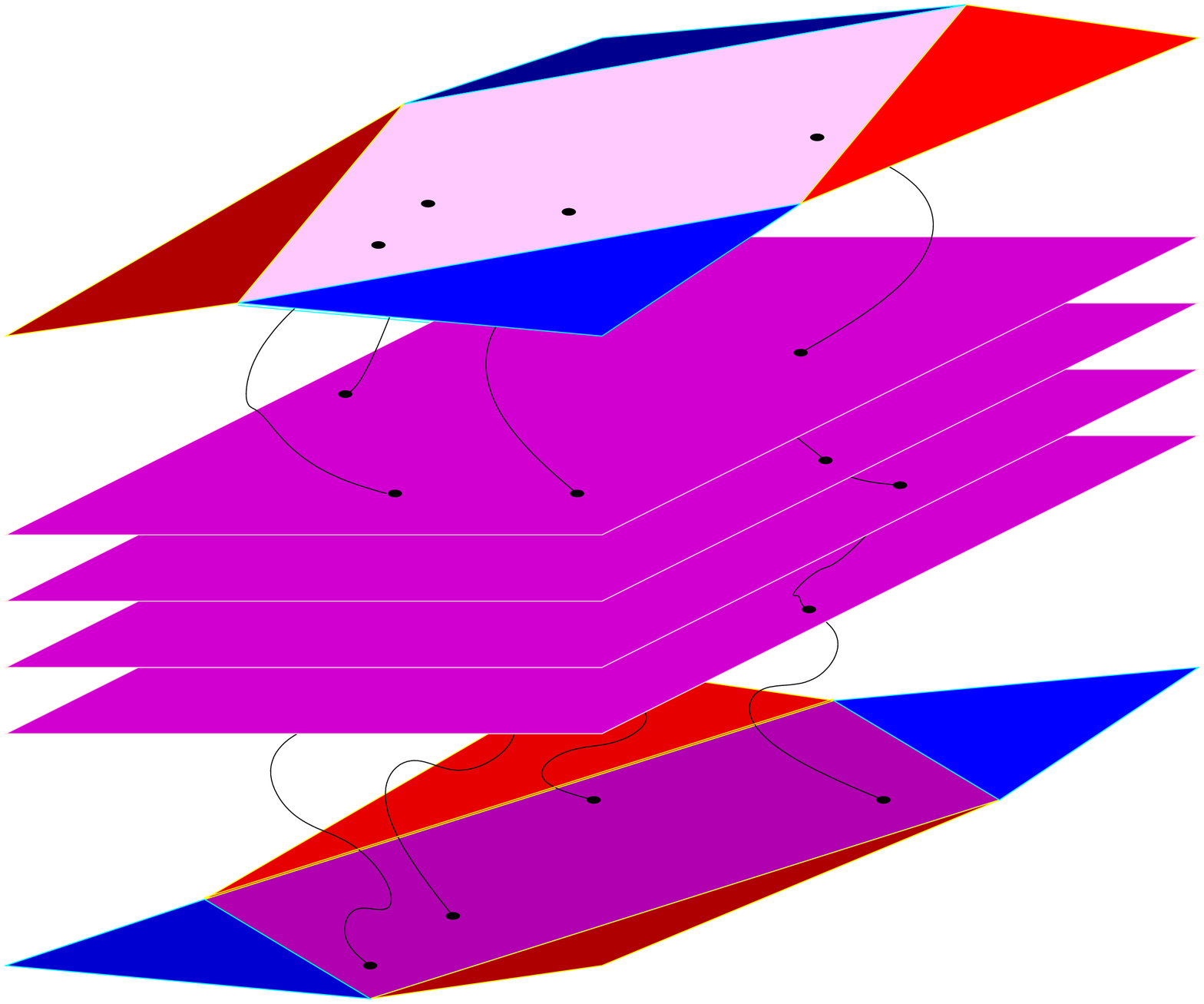,height=3in}}}
\put(218,110){\begin{rotate}{-20}{\Large $X_P$}\end{rotate}}
\put(300,80){{\LARGE \scalebox{1.8}[1]
{$\pmb{\twoheadrightarrow}$}}}
\put(328,30){$t\!=\!0$}
\put(340,42){\line(0,1){142}}
\put(340,42){\circle*{3}}
\put(340,184){\circle*{3}}
\put(344,96){\scalebox{1}[1.3]{{\LARGE $\Pro^1_\C$}}}
\put(328,186){$t\!=\!\infty$}
\end{picture} 

\bigskip
\bigskip
\bigskip
In particular, the only portion of the lower hull of $\tilde{P}$ 
(i.e., the ``lower portion'' of toric infinity on $X_{\tilde{P}}$) 
which is touched by zero set of $\hat{F}$ in $X_{\tilde{P}}$ is the 
parallelogram facet, and the projection of this facet has area 
exactly $\alpha\delta+\beta\gamma$. \dia 
\end{ex} 

\noindent
{\bf Proof of Theorem \ref{thm:hardbernie}:} We will first prove 
the generic case, and then derive the degenerate case, just 
as we did for the unmixed case.

At this point, we could just use Theorem \ref{thm:bernie} to get 
the generic case and proceed with our proof of the degenerate case. 
However, let us observe that we could instead use mixed subdivisions 
to directly obtain Theorem \ref{thm:bernie} without reducing to 
the unmixed case. The proof proceeds exactly like the proof of 
Theorem \ref{thm:kush}, except for the following differences: 
\begin{enumerate} 
\item{We work with $X_{\tilde{P}}$ instead of $Y_{\tilde{A}}$, 
where $\tilde{P}\!=\!(\hat{P}_1+\cdots+\hat{P}_n)\times [0,1]$.} 
\item{The map $\pi$ is essentially the same but is instead defined 
via the Cox coordinate ring \cite{cox}.} 
\item{The only portions of toric infinity in $X_{\tilde{P}}$ that 
intersect $\pi^{-1}(0)\cap \tilde{Z}$ are those corresponding 
to facets on the lower hull of $\tilde{P}$ that project to {\bf mixed cells} of 
$(P_1,\ldots,P_n)_{\omega}$. } 
\item{The final count of roots becomes a sum of roots of a collection of 
{\bf binomial} systems.} 
\end{enumerate} 

To prove the degenerate case, 
we proceed exactly as in the proof of Theorem \ref{thm:hardkush}, except 
with the following minor modifications: 
\begin{enumerate} 
\item{We use the notational changes above.}  
\item{The space of $F$ we work with is instead $\Pro^{N_1-1}_\C\times \cdots \times 
\Pro^{N_n-1}_\C$, where $N_i\!=\!\#A_i$ for all $i$. \qed}
\end{enumerate} 

We are now finally ready to state and prove the full version of Theorem 1: 

{\sc Theorem 1} (Full Version) 
{\em 
Suppose $F$ is any $k\times n$ polynomial system with 
$\supp(f_i)\!\subseteq\!A_i\!\subset\!(\N\cup\{0\})^n$ for all 
$i$ and let $Z_\C(F)$ denote the zero set of $F$ in $\Cn$. Then the 
number of connected components of $Z_\C(F)$ is no more than 
$\vol\left(\{\bO,e_1,\ldots,e_n)\cup\bigcup^k_{i=1} A_i\right)$ or 
$\cM\left(\{\bO,e_1\}\cup A_1,\ldots,\{\bO,e_k\}\cup A_k\right)$,  
according as $k\!<\!n$ or $k\!\geq\!n$. } 

\noindent
{\bf Proof:} Let $B\!:=\!\{\bO,e_1,\ldots,e_n)\cup\bigcup^k_{i=1} A_i$.  
If $k\!<\!n$ then we can simply set $f_{k+1}\!=\cdots =\!f_n\!=\!f_1$ and 
then apply Theorem \ref{thm:hardkush}, noting that the $\supp(f_i)\!\subseteq\!B$
for all $i$. In particular, it is easily checked that $Y_B$ actually 
contains an embedded copy of $\Cn$.  

To prove the case $k\!>\!n$, note that we can reconsider such an $F$ as a 
$k\times k$ polynomial system with $\supp(f_i)\!\subseteq\!\{0,e_i\}\cup A_i$ for 
all $i$. Once again, by virtue of the fact that 
$\bar{P}\!:=\!\conv(A_1)+\cdots+\conv(A_k)$ contains the non-negative orthant as 
one of the cones in its normal fan \cite{cox}, we have that our underlying 
toric variety (this time, $X_{\bar{P}}$) has an embedded copy of $\Cn$. 
So by Theorem \ref{thm:hardbernie}, we're done. \qed 

\bibliographystyle{acm}

\begin{thebibliography}{A}

\bibitem[Ber75]{bernie} Bernshtein, D. N., {\it ``The Number of
Roots of a System of Equations,"} Functional Analysis and its Applications
(translated from Russian), Vol. 9, No. 2, (1975),
pp.\ 183--185.

\bibitem[BCSS98]{bcss} Blum, Lenore; Cucker, Felipe; Shub, Mike; and 
Smale, Steve, {\it Complexity and Real Computation,} Springer-Verlag, 1998.         
\bibitem[Can88]{pspace} Canny, John F., {\it ``Some Algebraic
and Geometric Computations in PSPACE,''} Proc.\ 20$\thth$ ACM
Symp.\ Theory of Computing, Chicago (1988), ACM Press.

\bibitem[BCS97]{bcs} B\"urgisser, Peter; Clausen, Michael; and Shokrollahi, M.\ Amin, 
{\it Algebraic complexity theory,} with the collaboration of Thomas Lickteig, 
Grundlehren der Mathematischen Wissenschaften [Fundamental Principles of Mathematical 
Sciences], 315, Springer-Verlag, Berlin, 1997. 

\bibitem[Can93]{can93} \underline{\hspace{\jfc}}, {\it ``Computing roadmaps of 
general semi-algebraic sets,''} Comput.\ J.\ {\bf 36} (1993), no.\ 5, pp.\ 
504--514.

\bibitem[Cox03]{cox} Cox, David A., {\it 
``What is a Toric Variety,''} Tutorial for a conference on Algebraic Geometry and 
Geometric Modelling (Vilnius, Lithuania, July 29 -- August 2, 2002), submitted for 
publication, downloadable from {\tt http://www.cs.amherst.edu/~dac/lectures/tutorial.ps}  .  

\bibitem[CLO92]{clo1} Cox, David A., Little, John, and O'Shea, Donal, {\it 
Ideals, Varieties, and Algorithms,} Undergraduate Texts in Mathematics,
Springer-Verlag, 1992.

\bibitem[CLO98]{clo2} \underline{\hspace{\clo}}, {\it Using Algebraic
Geometry,} Graduate Texts in Mathematics 185, Springer-Verlag (1998).

\bibitem[Csa76]{csanky} Csanky, L., {\it ``Fast Parallel Matrix Inversion 
Algorithms,''}  SIAM J.\ Comput.\ {\bf 5} (1976), no.\ 4, pp.\ 618--623. 

\bibitem[EGA1]{ega1} Dieudonn\'e, Jean and Grothendieck, Alexander, {\it  
\'El\'ements de g\'eom\'etrie alg\'ebrique I: Le langage des 
sch\'emas,} Inst.\ Hautes \'Etudes Sci.\ Publ.\ Math.\ No.\ 
4, 1960. 

\bibitem[EGA2]{ega2} \underline{\hspace{\gd}}, {\it  
\'El\'ements de g\'eom\'etrie alg\'ebrique II: \'Etude globale 
\'el\'ementaire de quelques classes de morphismes,} 
Inst.\ Hautes \'Etudes Sci.\ Publ.\ Math.\ No.\ 
8, 1961. 

\bibitem[EGA3.1]{ega31} \underline{\hspace{\gd}}, {\it  
\'El\'ements de g\'eom\'etrie alg\'ebrique III: \'Etude 
cohomologique des faisceaux coh\'erents I,} 
Inst.\ Hautes \'Etudes Sci.\ Publ.\ Math.\ No.\ 
11, 1961. 

\bibitem[EGA3.2]{ega32} \underline{\hspace{\gd}}, {\it  
\'El\'ements de g\'eom\'etrie alg\'ebrique III: \'Etude 
cohomologique des faisceaux coh\'erents II,} 
Inst.\ Hautes \'Etudes Sci.\ Publ.\ Math.\ No.\ 
17, 1963. 

\bibitem[EGA4.1]{ega41} \underline{\hspace{\gd}}, {\it  
\'El\'ements de g\'eom\'etrie alg\'ebrique IV: \'Etude 
locale des sch\'emas et des morphismes de sch\'emas I,} 
Inst.\ Hautes \'Etudes Sci.\ Publ.\ Math.\ No.\ 
20, 1964. 

\bibitem[EGA4.2]{ega42} \underline{\hspace{\gd}}, {\it  
\'El\'ements de g\'eom\'etrie alg\'ebrique IV: \'Etude 
locale des sch\'emas et des morphismes de sch\'emas II,} 
Inst.\ Hautes \'Etudes Sci.\ Publ.\ Math.\ No.\ 
24, 1965. 

\bibitem[EGA4.3]{ega43} \underline{\hspace{\gd}}, {\it  
\'El\'ements de g\'eom\'etrie alg\'ebrique IV: \'Etude 
locale des sch\'emas et des morphismes de sch\'emas III,} 
Inst.\ Hautes \'Etudes Sci.\ Publ.\ Math.\ No.\ 
28, 1966. 

\bibitem[Dod01]{dod01} Dodis, Yevgeny, Personal Communication, e-mailed 
{}from Courant Institute, New York.

\bibitem[DGH98]{dgh} Dyer, Martin; Gritzmann, Peter; and
Hufnagel, Alexander, {\it ``On the Complexity of Computing Mixed Volumes,''}
SIAM J.\ Comput.\ {\bf 27} (1998), no.\ 2, pp.\ 356--400.

\bibitem[DMcK72]{dym} Dym, H.\ and McKean, H.\ P., {\it 
Fourier Series and Integrals,} Probability and Mathematical 
Statistics, vol.\ 14, Academic Press, 1972. 

\bibitem[EM99]{em99} Emiris, Ioannis Z.\ and Mourrain, Bernard, {\it 
``Computer algebra methods for studying and computing molecular conformations,''}
Algorithmica {\bf 25} (1999), no.\ 2--3, pp.\ 372--402.

\bibitem[EP02]{emirispan} Emiris, Ioannis Z.\ and Pan, Victor Y., {\it 
``Symbolic and Numeric Methods for Exploiting Structure in Constructing Resultant Matrices,''} 
 Journal of Symbolic Computation, Vol.\ 33, No.\ 4, April 1, 2002. 

\bibitem[FH95]{fh95}  Forsythe, Keith and Hatke, Gary, {\it ``A Polynomial 
Rooting 
Algorithm for Direction Finding,''} preprint, MIT Lincoln Laboratories, 1995.

\bibitem[Ful93]{tfulton} Fulton, William, {\it
Introduction to Toric Varieties}, Annals of Mathematics Studies, no.\ 131,
Princeton University Press, Princeton, New Jersey, 1993.

\bibitem[Gat01]{gat01} Gatermann, Karin, {\it ``Counting Stable Solutions of
Sparse Polynomial Systems in Chemistry,''}
Contemporary Mathematics, vol.\ 286, AMS-IMS-SIAM Joint Summer Research
Conference Proceedings of ``Symbolic Computation: Solving Equations in
Algebra, Geometry, and Engineering (June 11--15, 2000, Mount Holyoke
College),'' edited by R.\ Laubenbacher and V. Powers, pp.\ 53--69, AMS
Press, 2001.

\bibitem[GH99]{gh99} Gatermann, Karin and Huber, Birkett, {\it ``A Family of
Sparse Polynomial Systems Arising in Chemical Reaction
Systems,''} preprint SC 99-27 (August 1999), Konrad-Zuse-Zentrum
f\"{u}r Informationstechnik Berlin.

\bibitem[GKZ94]{gkz} Gel'fand, I. M., Kapranov, M. M., and
Zelevinsky, A. V., {\it Discriminants, Resultants and Multidimensional
Determinants,} Birkh\"auser, Boston, 1994.

\bibitem[Gol03]{goldman} Goldman, Ron, {\it ``Polar Forms in Geometric 
Modelling and Algebraic Geometry,''} presentation at a conference on Algebraic
Geometry and Geometric Modelling (Vilnius, Lithuania, July 29 -- August 2, 2002),
submitted for publication. 

\bibitem[GKP94]{gkp} Graham, R. L., Knuth, D. E., and Patashnik, O., {\it
Concrete Mathematics: A Foundation for Computer Science}, 2$^\nd$ edition,
Addison-Wesley, 1994.

\bibitem[GH94]{gh} Griffiths, Phillip and Harris, Joseph, 
{\it Principles of Algebraic Geometry,} 
Reprint of the 1978 original, Wiley Classics Library, 
John Wiley \& Sons, Inc., New York, 1994. 

\bibitem[SGA1]{sga1} Grothendieck, Alexander, et.\ al., {\it  
Rev\^{e}tements \'etales et groupe fondamental,} 
S\'eminaire de G\'eom\'etrie Alg\'ebrique du Bois-Marie, 1960--1961, directed 
by Alexandre Grothendieck, with essays by M.\ Raynaud, Lecture Notes in 
Mathematics, Vol.\ 224, Springer-Verlag, Berlin-New York, 1971. 

\bibitem[SGA2]{sga2} \underline{\hspace{\groth}}, {\it Cohomologie locale des 
faisceaux coh\'erents et th\'eor\`emes de Lefschetz locaux et globaux 
$(SGA$ $2)$,} with an essay by Mich\`ele Raynaud, S\'eminaire de 
G\'eom\'etrie Alg\'ebrique du Bois-Marie, 1962, Advanced Studies in Pure 
Mathematics, Vol.\ 2, North-Holland Publishing Co., Amsterdam; Masson \& Cie, 
\'Editeur, Paris, 1968. 

\bibitem[SGA3.1]{sga31} \underline{\hspace{\groth}}, {\it  
Sch\'emas en groupes, vol.\ I: Propri\'et\'es g\'en\'erales des sch\'emas en 
groupes,}  S\'eminaire de G\'eom\'etrie Alg\'ebrique du Bois-Marie, 1962/1964, 
directed by M.\ Demazure and A.\ Grothendieck, Lecture Notes in Mathematics, 
Vol.\ 151, Springer-Verlag, Berlin-New York 1962/1964. 

\bibitem[SGA3.2]{sga32} \underline{\hspace{\groth}}, {\it  
Sch\'emas en groupes, vol.\ II: Groupes de type multiplicatif, et structure 
des sch\'emas en groupes g\'en\'eraux,}  
S\'eminaire de G\'eom\'etrie Alg\'ebrique du Bois-Marie, 1962/1964, directed 
by M.\ Demazure and A.\ Grothendieck, Lecture Notes in Mathematics, Vol.\ 152, 
Springer-Verlag, Berlin-New York 1962/1964. 

\bibitem[SGA3.3]{sga33} \underline{\hspace{\groth}}, {\it  
Sch\'emas en groupes, vol.\ III: Structure des sch\'emas en groupes 
r\'eductifs,}  
S\'eminaire de G\'eom\'etrie Alg\'ebrique du Bois-Marie, 1962/1964, directed 
by M.\ Demazure and A.\ Grothendieck, Lecture Notes in Mathematics, Vol.\ 153, 
Springer-Verlag, Berlin-New York 1962/1964. 

\bibitem[SGA4.1]{sga41} \underline{\hspace{\groth}}, {\it 
Th\'eorie des topos et cohomologie \'etale des sch\'emas, vol.\ 1: Th\'eorie 
des topos} 
S\'eminaire de G\'eom\'etrie Alg\'ebrique du Bois-Marie, 1963--1964, directed 
by M.\ Artin, A.\ Grothendieck, and J.\ L.\ Verdier, with the collaboration of 
N.\ Bourbaki, P.\ Deligne and B.\ Saint-Donat, Lecture Notes in Mathematics, 
Vol.\ 269, Springer-Verlag, Berlin-New York, 1972.

\bibitem[SGA4.2]{sga42} \underline{\hspace{\groth}}, {\it 
Th\'eorie des topos et cohomologie \'etale des sch\'emas, vol.\ 2} 
S\'eminaire de G\'eom\'etrie Alg\'ebrique du Bois-Marie, 1963--1964, directed 
by M.\ Artin, A.\ Grothendieck, and J.\ L.\ Verdier, with the collaboration of 
N.\ Bourbaki, P.\ Deligne and B.\ Saint-Donat, Lecture Notes in Mathematics, 
Vol.\ 270, Springer-Verlag, Berlin-New York, 1972.

\bibitem[SGA4.3]{sga43} \underline{\hspace{\groth}}, {\it 
Th\'eorie des topos et cohomologie \'etale des sch\'emas, vol.\ 3} 
S\'eminaire de G\'eom\'etrie Alg\'ebrique du Bois-Marie, 1963--1964, directed 
by M.\ Artin, A.\ Grothendieck, and J.\ L.\ Verdier, with the collaboration of 
P.\ Deligne and B.\ Saint-Donat, Lecture Notes in Mathematics, Vol.\ 305, 
Springer-Verlag, Berlin-New York, 1973.

\bibitem[SGA4$\frac{1}{2}$]{sga45} Deligne, 
Pierre, {\it Cohomologie \'etale,} 
S\'eminaire de G\'eom\'etrie Alg\'ebrique du Bois-Marie SGA 4\scalebox{1}[.7]
{$\frac{1}{2}$}, with the 
collaboration of J.\ F.\ Boutot, A.\ Grothendieck, L.\ Illusie and J.\ L.\ 
Verdier, Lecture Notes in Mathematics, Vol.\ 569, Springer-Verlag, Berlin-New 
York, 1977.

\bibitem[SGA5]{sga5} Grothendieck, Alexander, et.\ al., {\it 
Cohomologie $l$-adique et fonctions $L$,} unpublished. 

\bibitem[SGA6]{sga6} \underline{\hspace{\groth}}, {\it 
Th\'eorie des intersections et th\'eor\`eme de Riemann-Roch,}  
S\'eminaire de G\'eom\'etrie Alg\'ebrique du Bois-Marie, 1966--1967, directed 
by P.\ Berthelot, A.\ Grothendieck, and L.\ Illusie, with the collaboration of 
D.\ Ferrand, J.\ P.\ Jouanolou, O.\ Jussila, S.\ Kleiman, M.\ Raynaud and 
J.\ P.\ Serre, Lecture Notes in Mathematics, Vol.\ 225, 
Springer-Verlag, Berlin-New York, 1971.

\bibitem[SGA7]{sga7} \underline{\hspace{\groth}},  {\it Groupes de monodromie 
en g\'eom\'etrie alg\'ebrique I,}  
S\'eminaire de G\'eom\'etrie Alg\'ebrique du Bois-Marie, 1967--1969, directed 
by A.\ Grothendieck, with the collaboration of M.\ Raynaud and D.\ S.\ Rim, 
Lecture Notes in Mathematics, Vol.\ 288, Springer-Verlag, Berlin-New York, 
1972.

\bibitem[HMPS00]{hmps} H\"agele, Klemens; Morais, Juan Enrique; Pardo, Luis
Miguel; Sombra, Martin, {\it ``On the Intrinsic Complexity of the
Arithmetic Nullstellensatz,''} Journal of Pure and Applied Algebra
{\bf 146} (2000), no.\ 2, pp.\ 103--183.

\bibitem[Har77]{hartshorne} Hartshorne, Robin, {\it Algebraic
Geometry,''} Graduate Texts in Mathematics, No.\ 52,
Springer-Verlag.

\bibitem[Hir94]{hirsch} Hirsch, Morris, {\it Differential Topology,}
corrected reprint of the 1976 original, Graduate Texts in Mathematics,
33, Springer-Verlag, New York, 1994.

\bibitem[HS95]{hs} Huber, Birk and Sturmfels, Bernd, {\it ``A Polyhedral 
Method for Solving Sparse Polynomial Systems,''} Math.\ Comp.\ {\bf 64} (1995),
no.\ 212, pp.\ 1541--1555.

\bibitem[Ili89]{unimod} Iliopoulos, Costas S., {\it ``Worst Case
Complexity Bounds on Algorithms for Computing the Canonical Structure of
Finite Abelian Groups and the Hermite and Smith Normal Forms of an Integer
Matrix,"} SIAM Journal on Computing, 18 (1989), no.\ 4, pp.\ 658--669.

\bibitem[JKSS03]{jkss} Jeronimo, Gabriela; Krick, Teresa; Sabia, Juan; and 
Sombra, Martin, {\it ``The Computational Complexity of the Chow Form,''} 
Math ArXiV preprint {\tt math.AG/0210009}. 

\bibitem[KLS97]{kls} Kannan, Ravi; Lov\'asz, L\'aszl\'o; and Simonovitz, Miki,
{\it ``Random Walks and an $\cO^*(n^5)$ Volume Algorithm
for Convex Bodies,"}
Random Structures Algorithms, {\bf 11} (1997), no.\ 1, pp.\ 1--50.

\bibitem[KM97]{km97a} Karpinski, Marek and Macintyre, Angus J., {\it 
``Polynomial
bounds for VC dimension of sigmoidal and general Pfaffian neural networks,''}
{\em J.\ Comp.\ Sys.\ Sci.}, 54, pp.\ 169--176, 1997.

\bibitem[Khe02]{khetan} Khetan, Amit, {\it ``Determinental Formula 
for the Chow Form of a Toric Surface'',} Proceedings of 
the International Symposium on Symbolic and Algebraic Computation 
(ISSAC) 2002, ACM Press, to appear.  

\bibitem[Kho91]{few} \underline{\hspace{\khov}}, {\it Fewnomials,}
AMS Press, Providence, Rhode Island, 1991.

\bibitem[Koi96]{hnam} Koiran, Pascal, {\it ``Hilbert's Nullstellensatz
is in the Polynomial Hierarchy,''} DIMACS Technical Report 96-27,
July 1996. ({\bf Note:} This preprint considerably improves the published
version which appeared in Journal of Complexity in 1996.)

\bibitem[Koi97]{koiran} \underline{\hspace{\koi}}, {\it ``Randomized and
Deterministic Algorithms for the Dimension of Algebraic Varieties,''}
Proceedings of the 38$\thth$ Annual IEEE Computer Society
Conference on Foundations of Computer Science (FOCS),
Oct.\ 20--22, 1997, ACM Press.

\bibitem[KPS01]{cool} Krick, Teresa; Pardo, Luis Miguel; and Sombra, Martin,
{\it ``Sharp Arithmetic Nullstellensatz,''} Duke Mathematical Journal 
{\bf 109} (2001), no.\ 3, pp.\ 521--598. 

\bibitem[Kus75]{kush1} Kushnirenko, Anatoly Georievich, {\it ``A Newton 
Polytope and the Number of Solutions of a System of $k$ Equations in $k$ 
Unknowns,"} Usp.\ Matem.\ Nauk., 30, no.\ 2, pp.\ 266--267 (1975).

\bibitem[Kus76]{kush2} \underline{\hspace{\agku}}, {\it ``Newton Polytopes and
the B\'ezout Theorem,"} Functional Analysis and its Applications (translated
{}from Russian), vol.\ 10, no.\ 3, July--September (1976), pp.\ 82--83.

\bibitem[Lec00]{lecerf} Lecerf, Gr\'egoire, {\it ``Computing an 
Equidimensional Decomposition of an Algebraic Variety by Means 
of Geometric Resolutions,''} Proceedings of the 
International Symposium on Symbolic Algebra and Computation 
(ISSAC) 2000.

\bibitem[Li97]{li97}  Li, Tien Yien, {\it ``Numerical solution of multivariate polynomial 
systems by homotopy continuation methods,''} Acta numerica, 1997, pp.\ 399--436, Acta Numer., 
6, Cambridge Univ.\ Press, Cambridge, 1997. 

\bibitem[LRW03]{tri} Li, Tien-Yien; Rojas, J.\ Maurice; and
Wang, Xiaoshen, {\it ``Counting Real Connected Components of Trinomial
Curves Intersections and m-nomial Hypersurfaces,''} Discrete and
Computational Geometry, to appear.

\bibitem[LW91]{lw91} Li, Tien Yien and Wang, Xiaoshen, {\it `` Solving deficient polynomial 
systems with homotopies which keep the subschemes at infinity invariant,''} 
Math.\ Comp.\ 56 (1991), no.\ 194, pp.\ 693--710.

\bibitem[MR03]{high} Malajovich, Gregorio and Rojas, J.\ Maurice, 
{\it ``High 
Probability Analysis of the Condition Number of Sparse Polynomial Systems,''} 
submitted for publication, also available as math ArXiV preprint 
{\tt math.NA/0212179}.

\bibitem[Man98]{man98} Manocha, Dinesh, {\it ``Numerical Methods for Solving 
Polynomial
Equations,''} Applications of Computational Algebraic Geometry (San Diego, CA,
1997), pp.\ 41--66, Proc.\ Sympos.\ Appl.\ Math., 53,
Amer.\ Math.\ Soc., Providence, RI, 1998.

\bibitem[McD02]{mcdonald} McDonald, John, {\it ``Fractional power series solutions for 
systems of equations,''} 
Discrete Comput.\ Geom.\ 27 (2002), no.\ 4, pp.\ 501--529.

\bibitem[McL97]{mcl97} McLennan, Andrew, {\it ``The maximal number of regular 
totally mixed Nash equilibria,''} J.\ Econom.\ Theory 72 (1997), no.\ 2, pp.\ 
411--425.

\bibitem[MS87]{ms} Morgan, Alexander and Sommese, Andrew, 
{\it ``A homotopy for solving general polynomial systems that respects $m$-homogeneous 
structures,''} Appl.\ Math.\ Comput.\ 24 (1987), no.\ 2, pp.\ 101--113.

\bibitem[Mou02]{mourrain} Mourrain, Bernard, 
{\it ``Resultants,''} presentation at a conference on Algebraic
Geometry and Geometric Modelling (Vilnius, Lithuania, July 29 -- August 2, 2002),
submitted for publication. 

\bibitem[MP98]{mourrainpan} Mourrain, Bernard and Pan, Victor,
{\it ``Asymptotic Acceleration of Solving Multivariate Polynomial
Systems of Equations,''} Proc.\ STOC '98, pp.\ 488--496, ACM Press, 1998.

\bibitem[Mum95]{mumford1} Mumford, David, {\it Algebraic Geometry I: Complex 
Projective Varieties,} reprint of the 1976 edition, Classics in Mathematics,  
Springer-Verlag, Berlin, 1995.

\bibitem[Mum99]{mumford2} \underline{\hspace{\mum}}, {\it  The red book of 
varieties and schemes,} second, expanded edition, includes the Michigan 
lectures (1974) on curves and their Jacobians, with contributions by Enrico 
Arbarello, Lecture Notes in Mathematics, 1358, Springer-Verlag, Berlin, 1999.

\bibitem[NR96]{nr96} Neff, C.\ Andrew and Reif, John, {\it ``An
Efficient Algorithm for the Complex Roots Problem,''}
Journal of Complexity {\bf 12} (1996), no.\ 2, pp.\ 81--115.

\bibitem[NM99]{nm99} Ne\v{s}i\'c, D.\ and Mareels, Ivan M.\ Y.,
{\it ``Controllability of structured polynomial systems,''}
IEEE Trans. Automat. Control 44 (1999), no. 4, pp.\ 761--764.

\bibitem[Pap95]{papa} Papadimitriou, Christos H., {\it Computational
Complexity,} Addison-Wesley, 1995.

\bibitem[Pla84]{plaisted} Plaisted, David A., {\it ``New NP-Hard and
NP-Complete Polynomial and Integer Divisibility Problems,''}
Theoret.\ Comput.\ Sci.\ 31 (1984), no.\ 1--2, 125--138.

\bibitem[PS85]{preparata} Preparata, Franco P. and
Shamos, Michael Ian, {\it Computational Geometry: An Introduction,}
Texts and Monographs in Computer Science, Springer-Verlag,
New York-Berlin, 1985.

\bibitem[Roj94]{convex} Rojas, J.\ Maurice, {\it ``A Convex Geometric
Approach to Counting the Roots of a Polynomial
System,"} Theoretical Computer Science (1994), vol.\ 133 (1), pp.\ 105--140.
(Additional notes and corrections available on-line at {\tt
http://www.math.tamu.edu/\~rojas/list2.html} .)

\bibitem[Roj97]{roj97} \underline{\hspace{\jmr}}, {\it
``A New Approach to Counting Nash Equilibria,''} Proceedings of the IEEE/IAFE
Conference on Computational Intelligence for Financial Engineering, Manhattan,
New York, March 23--25, 1997, pp.\ 130--136.

\bibitem[Roj99a]{toric} \underline{\hspace{\jmr}}, {\it ``Toric
Intersection Theory for Affine Root Counting,''} Journal of Pure and
Applied Algebra, vol.\ 136, no.\ 1, March, 1999, pp.\ 67--100.

\bibitem[Roj99b]{gcp} \underline{\hspace{\jmr}}, {\it ``Solving Degenerate
Sparse Polynomial Systems Faster,''} Journal of Symbolic Computation,
vol.\ 28 (special issue on
elimination theory), no.\ 1/2, July and August 1999, pp.\ 155--186.

\bibitem[Roj00a]{four} \underline{\hspace{\jmr}}, {\it ``Algebraic 
Geometry Over Four Rings and the Frontier to Tractability,''} Contemporary 
Mathematics, vol.\ 270, Proceedings of a Conference on Hilbert's Tenth Problem 
and Related Subjects (University of Gent, November 1-5, 1999), edited by
Jan Denef, Leonard Lipschitz, Thanases Pheidas, and Jan Van
Geel, pp.\ 275--321, AMS Press (2000).                     

\bibitem[Roj00b]{real}  \underline{\hspace{\jmr}}, ``Some Speed-Ups and Speed 
Limits for Real Algebraic Geometry,'' Journal of Complexity, FoCM 1999 
special issue, vol.\ 16, no.\ 3 (sept.\ 2000), pp.\ 552--571. 

\bibitem[Roj01]{jcs} \underline{\hspace{\jmr}}, {\it ``Computational
Arithmetic Geometry I: Sentences Nearly in the Polynomial Hierarchy,''}
J.\ Comput.\ System Sci., STOC '99 special issue, vol.\ 62, no.\ 2,
march 2001, pp.\ 216--235. 

\bibitem[Roj02]{add} \underline{\hspace{\jmr}}, {\it ``Additive
Complexity and the Roots of Polynomials Over
Number Fields and $\cp$-adic Fields,''} Proceedings of
the 5$\thth$ Annual Algorithmic Number Theory Symposium
(ANTS V), Lecture Notes in Computer Science \#2369, pp.\ 506--515,
Springer-Verlag (2002).

\bibitem[Roj03]{ari} \underline{\hspace{\jmr}}, {\it ``Arithmetic
Multivariate Descartes' Rule,''} Math ArXiV preprint {\tt math.NT/0110327},
submitted for publication. 

\bibitem[RY02]{ry} Rojas, J.\ M.\ and Ye, Yinyu, {\it ``On Solving
Fewnomials Over an Interval in Fewnomial Time,''} submitted for publication,
also available as Math ArXiV preprint {\tt math.NA/0106225}.

\bibitem[Sha94]{sha} Shafarevich, Igor R., {\it Basic
Algebraic Geometry I,} second edition, Springer-Verlag (1994).

\bibitem[Shu93]{shub} Shub, Mike, {\it ``Some Remarks
on B\'ezout's Theorem and Complexity Theory,''} {}From
Topology to Computation: Proceedings of
the Smalefest (Berkeley, 1990), pp.\ 443--455, Springer-Verlag, 1993.   

\bibitem[Sma00]{smale} Smale, Steve,  {\it ``Mathematical
Problems for the Next Century,''} Mathematics: 
Frontiers and Perspectives, pp.\ 271--294, 
Amer.\ Math.\ Soc., Providence, RI, 2000. 

\bibitem[Smi61]{smith} Smith, H.\ J.\ S., {\it ``On Systems of Integer
Equations and Congruences,''} Philos.\ Trans.\ 151, pp.\ 293--326 (1861).

\bibitem[Sot03]{frank} Sottile, Frank, {\it ``Toric Ideals, 
Real Toric Varieties, and the Moment Map,''} presentation at a conference on Algebraic 
Geometry and Geometric Modelling (Vilnius, Lithuania, July 29 -- August 2, 2002),
submitted for publication, also available as Math ArXiV 
preprint {\tt math.AG/0212044}. 

\bibitem[Str98]{strang} Strang, Gilbert, {\it Introduction to Linear Algebra,} 
Wellesley-Cambridge Press, 1998. 

\bibitem[Sus98]{sus98} Sussmann, H\'ector J., {\it ``Some Optimal Control 
Applications
of Real-Analytic Stratifications and Desingularization,''}
Singularities Symposium --- \L{}ojasiewicz 70 (Krak\'ow, 1996; Warsaw, 1996),
211--232, Banach Center Publ., 44,

\bibitem[Ver00]{verschelde} Verschelde, Jan, {\it ``Toric Newton method for polynomial 
homotopies,''} Symbolic computation in algebra, analysis, and geometry (Berkeley, CA, 1998), 
J.\ Symbolic Comput.\ 29 (2000), no.\ 4--5, pp.\ 777--793.

\bibitem[Vid97]{vid97} Vidyasagar, M., {\it
A theory of learning and generalization,}
With applications to neural networks and control systems. Communications and
Control Engineering Series, Springer-Verlag London, Ltd., London, 1997.

\bibitem[VR02]{vr02} Vidyasagar, M.\ and Rojas, J.\ Maurice, {\it ``An 
Improved Bound 
on the VC-Dimension of Neural Networks with Polynomial Activation Functions,''}
submitted for publication.

\end{thebibliography}

\end{document}